\documentclass[12pt]{amsart}
\usepackage{amssymb}
\usepackage{amsfonts}
\usepackage{amscd}

\swapnumbers

\newtheorem{theorem}[subsubsection]{Theorem}
\newtheorem{lem}[subsubsection]{Lemma}
\newtheorem{cor}[subsubsection]{Corollary}
\newtheorem{prop}[subsubsection]{Proposition}

\newtheorem{itheorem}[subsection]{Theorem}

\theoremstyle{definition}

\newtheorem{def-theorem}[subsubsection]{Definition-Theorem}
\newtheorem{def-prop}[subsubsection]{Definition-Proposition}
\theoremstyle{remark}
\newtheorem{remark}[subsubsection]{Remark}

\makeatletter
\def\notocsection{\@startsection{section}{\@M}%
  \z@{1.5\linespacing\@plus\linespacing}{.5\linespacing}%
  {\normalfont\bfseries\centering}}
\makeatother

\theoremstyle{plain}
\numberwithin{equation}{section}

\def\boxit#1#2{\setbox1=\hbox{\kern#1{#2}\kern#1}%
\dimen1=\ht1 \advance\dimen1 by #1
\dimen2=\dp1 \advance\dimen2 by #1
\setbox1=\hbox{\vrule height\dimen1 depth\dimen2\box1\vrule}%
\setbox1=\vbox{\hrule\box1\hrule}%
\advance\dimen1 by .4pt \ht1=\dimen1
\advance\dimen2 by .4pt \dp1=\dimen2 \box1\relax}

\let\cal\mathcal

\def\AA{{\mathbf A}}
\def\BB{{\mathbf B}}
\def\CC{{\mathbf C}}

\def\FF{{\mathbf F}}
\def\GG{{\mathbf G}}

\def\LL{{\mathbf L}}

\def\NN{{\mathbf N}}

\def\PP{{\mathbf P}}
\def\QQ{{\mathbf Q}}
\def\RR{{\mathbf R}}

\def\XX{{\mathbf X}}
\def\YY{{\mathbf Y}}
\def\ZZ{{\mathbf Z}}

\def\aa{{\mathbf a}}
\def\bb{{\mathbf b}}

\def\spec{{\rm Spec}\,}

\def\cA{{\mathcal A}}
\def\cB{{\mathcal B}}
\def\cC{{\mathcal C}}

\def\cF{{\mathcal F}}

\def\cI{{\mathcal I}}

\def\cL{{\mathcal L}}

\def\cO{{\mathcal O}}

\def\cW{{\mathcal W}}
\def\cX{{\mathcal X}}

\def\cZ{{\mathcal Z}}
\mathchardef\alphag="7C0B
\mathchardef\betag="7C0C
\mathchardef\gammag="7C0D
\mathchardef\deltag="7C0E
\mathchardef\varepsilong="7C22
\mathchardef\varphig="7C27
\mathchardef\psig="7C20
\mathchardef\zetag="7C10
\mathchardef\epsilong="7C0F
\mathchardef\rhog="7C1A
\mathchardef\taug="7C1C
\mathchardef\upsilong="7C1D
\mathchardef\iotag="7C13
\mathchardef\thetag="7C12
\mathchardef\pig="7C19
\mathchardef\sigmag="7C1B
\mathchardef\etag="7C11
\mathchardef\omegag="7C21
\mathchardef\kappag="7C14
\mathchardef\lambdag="7C15
\mathchardef\mug="7C16
\mathchardef\xig="7C18
\mathchardef\chig="7C1F
\mathchardef\nug="7C17
\mathchardef\varthetag="7C23
\mathchardef\varpig="7C24
\mathchardef\varrhog="7C25
\mathchardef\varsigmag="7C26
\mathchardef\Omegag="7C0A
\mathchardef\Thetag="7C02
\mathchardef\Sigmag="7C06
\mathchardef\Deltag="7C01
\mathchardef\Phig="7C08
\mathchardef\Gammag="7C00
\mathchardef\Psig="7C09
\mathchardef\Lambdag="7C03
\mathchardef\Xig="7C04
\mathchardef\Pig="7C05
\mathchardef\Upsilong="7C07

\DeclareMathOperator*{\Spec}{Spec}
\def\ord{{\rm ord}}
\DeclareMathOperator*{\ac}{{\overline {\rm ac}}}

\begin{document}

\title[Definable sets, motives and
$P$-adic integrals]{Definable sets, motives and 
$P$-adic integrals}

\author{Jan Denef}
\address{University of Leuven, Department of Mathematics,
Celestijnenlaan 200B, 3001 Leu\-ven, Bel\-gium }
\email{ Jan.Denef@wis.kuleuven.ac.be}
\urladdr{http://www.wis.kuleuven.ac.be/wis/algebra/denef.html}

\author{Fran\c cois Loeser}

\address{Centre de Math{\'e}matiques,
Ecole Polytechnique,
F-91128 Palaiseau
(UMR 7640 du CNRS), {\rm and}
Institut de Math{\'e}matiques,
Universit{\'e} P. et M. Curie, Case 82,
4 place Jussieu,
F-75252 Paris Cedex 05
(UMR 7596 du CNRS)}
\email{loeser@math.polytechnique.fr}
\urladdr{http://math.polytechnique.fr/cmat/loeser/loeser.html}

%    General info

%\date{Revised \today}

%\keywords{Algebraic geometry}

%\begin{abstract}We associate a canonical virtual motive
%to definable sets over a field of characteristic zero.
%We use this construction to to show that 
%very general $p$-adic integrals are canonically interpolated by
%motivic ones.
%\end{abstract}

\maketitle

\notocsection{Introduction}

\subsection{}Let $X$ be a scheme, reduced and separated, of finite type
over $\ZZ$.
For $p$ a prime number, one may consider the set $X (\ZZ_{p})$
of its $\ZZ_{p}$-rational
points. 
For every $n$ in $\NN$, there is a natural map
$\pi_{n} : X (\ZZ_{p}) \rightarrow X (\ZZ / p^{n + 1})$ assigning
to a $\ZZ_{p}$-rational
point its class modulo $p^{n + 1}$.
The image $Y_{n, p}$ of $X (\ZZ_{p})$ by $\pi_{n}$ is exactly the set of
$\ZZ / p^{n + 1}$-rational points which can be lifted to 
$\ZZ_{p}$-rational
points. Denote by 
$N_{n, p}$ the cardinality of the finite set $Y_{n, p}$.
By a result of the first author \cite{rat},
the Poincar{\'e} series
\begin{equation}
P_{p} (T) := \sum_{n \in \NN} \, N_{n, p} \, T^{n}
\end{equation}
is a rational function of $T$. Later
Macintyre
\cite{Macintyre}, Pas \cite{P}
and Denef \cite{AMJ} proved that
the degrees of the numerator and denominator of the rational
function $P_{p} (T)$
are bounded independently of $p$. One task of the present paper is to
prove a much stronger uniformity result by constructing
a canonical rational function $P_{{\rm ar}} (T)$ which
specializes to $P_{p} (T)$ for almost all $p$.
It follows in particular from our results, that there exist,
for every $n$ in $\NN$,
varieties $Z_{n, i}$ and rational numbers $r_{n, i}$ in $\QQ$,
$1 \leq i \leq m_{n}$, such that, for almost all $p$ and every $n$,
\begin{equation}
N_{n, p}  = \sum_{1 \leq i \leq m_{n}} r_{n, i} \, {\rm card}
\, Z_{n, i} (\FF_{p}).
\end{equation}

Hence a natural idea would be to try to construct $P_{{\rm ar}} (T)$
as a series with coefficients in $K_{0} ({\rm Sch}_{\QQ}) \otimes \QQ$,
with $K_{0} ({\rm Sch}_{\QQ})$
the ``Grothendieck ring
of algebraic varieties over $\QQ$'', defined in \ref{1.2}.
However, since different varieties
over a number field may have the same $L$-function, and
we want the function $P_{{\rm ar}} (T)$ to be canonical, we have to replace
the varieties $Z_{n, i}$ by Chow motives and 
the naive Grothendieck ring
$K_{0} ({\rm Sch}_{\QQ}) \otimes \QQ$ by 
the ring $K_{0}^{v} ({\rm Mot}_{\QQ, \bar \QQ}) \otimes \QQ$,
with $K_{0}^{v} ({\rm Mot}_{\QQ, \bar \QQ})$ the image of
$K_{0} ({\rm Sch}_{\QQ})$ in the 
Grothendieck ring of Chow motives with coefficients in $\bar \QQ$,
as defined in \ref{Chow}.

We can now state our uniformity result on the series $P_{p} (T)$
as follows.
\begin{itheorem}Given $X$ as above, there exists
a canonical series $P_{{\rm ar}} (T)$
with coefficients in 
$K^{v}_{0} ({\rm Mot}_{\QQ, \bar \QQ}) \otimes {\QQ}$, which is a 
rational function of $T$ and
which
specializes - after taking the trace of the Frobenius on an
{\'e}tale realization - onto the $p$-adic Poincar{\'e} series
$P_{p} (T)$, for almost all
prime numbers $p$.
\end{itheorem}
In fact, one constructs in \ref{ar}, for $X$ a variety over any field $k$
of characteristic zero, a canonical series $P_{{\rm ar}} (T)$
with coefficients in the Grothendieck ring
$K_{0}^{v} ({\rm Mot}_{k, \bar \QQ}) \otimes {\QQ}$ - 
defined in \ref{Chow} -
which is equal to the former one when $k  = \QQ$ and,
furthermore, the series $P_{{\rm ar}} (T)$ is rational in a precise sense (Theorem
\ref{arrat}).

\subsection{}Assume now $X$ is a subscheme 
of the affine space
$\AA^{m}_{\ZZ}$  defined by $f_{1} = \ldots = f_{r} = 0$,
with $f_{1}$, \dots, $f_{r}$ in 
$\ZZ [x_{1}, \ldots, x_{m}]$.
The starting point in the 
proof of the rationality of the series
$P_{p} (T)$ in \cite{rat} is
the relation
$P_{p} (p^{- m - 1}p^{-s}) = \frac{p}{p-1} I_{p} (s)$,
with 
$$
I_{p} (s) = \int_{W_{p}} \vert w\vert^{s}_{p} |dx|_{p}|dw|_{p},
$$
where $W_{p}$ is defined as
\begin{multline*}
W_{p} := \Bigl\{(x, w) \in \ZZ_{p}^{m} \times \ZZ_{p}
\Bigm| \\ \exists
y \in \ZZ_{p}^{m}:
x \equiv y \mod w, \, \text{and} \, f_{i} (y) = 0, \, \text{for} \, i = 1, \ldots r
\Bigr\}.
\end{multline*}

In particular, $W_{p}$ is defined by a formula in the first
order language of valued
fields independently of the prime $p$.

More generally, let
$k$ be a finite
extension of $\QQ$ with ring of integers $\cO$ and $R = \cO
[\frac{1}{N}]$,
with $N$ a non zero multiple of the discriminant.
For $x$ a closed point of $\spec R$,
we denote by $K_{x}$ the completion of the localization of $R$ at $x$,
by
$\cO_{x}$ its ring of integer, 
and by $\FF_{x}$ the residue field at $x$, a finite field of cardinality
$q_{x}$.
Let 
$f (x)$ be a polynomial in
$k [x_{1}, \ldots, x_{m}]$ (or more generally a
definable function in the first order
language of valued fields
with variables and values taking place in  the valued field
and
with coefficients in $k$)
and let 
$\varphi$
be a formula in the language of valued fields with 
coefficients in $k$, free 
variables $x_{1}, \ldots, x_{m}$ running
over the valued field and no other free
variables.
Now set
$W_{x} := \{y \in \cO_{x}^m \vert \varphi (y) \, \text{holds}\}$
and define 
$$
I_{\varphi, f, x} (s) = \int_{W_{x}} \vert f\vert^{s}_{x} |dx|_{x},
$$
for $x$ a closed point of ${\rm Spec} \, \cO$, with similar notation as before.
In Theorem \ref{red}
we show more generally, that
there exists a canonical
rational function
$I_{\varphi, f, {\rm mot}} (T)$ with coefficients
in a appropriate Grothendieck ring of motives,
such that,
for almost all closed points $x$ in $\spec \cO$,
$I_{\varphi, f, {\rm mot}} (T)$
specializes - after taking the trace of the Frobenius on an
{\'e}tale realization and after setting
$T = q_{x}^{-s}$ - to the $p$-adic integral
$I_{\varphi, f, x} (s)$.

\subsection{}Let $k$ be a field of characteristic zero. In
section \ref{assi}, to every formula
$\varphi$ in the first order language
of rings
with coefficients in $k$
we associate, in a canonical way,
a virtual motive $\chi_{c} (\varphi)$
in the ring
$K_{0}^{v} ({\rm Mot}_{k, \bar \QQ}) \otimes {\QQ}$. This virtual motive 
$\chi_{c} (\varphi)$ depends only
on the sets defined by $\varphi$ over the
pseudo-finite fields containing $k$. When $k$ is of finite type over
$\QQ$, one recovers the number of points of the sets defined by
$\varphi$ over the residual fields at almost
all finite places of $k$, by taking the trace of the corresponding
Frobenius on an
{\'e}tale realization of
$\chi_{c} (\varphi)$. It is quite amazing to note that, since virtual
motives have Euler characteristic and Hodge numbers (with compact
supports), it follows from our results that every 
formula
with coefficients in $k$ has an Euler characteristic in $\QQ$ and also
Hodge numbers in $\QQ$. One of  the simplest example is the formula 
$(\exists y) \, x = y^{n}$ which has Euler characteristic
$\frac{1}{n}$ and 
Hodge polynomial $\frac{1}{n}uv$.

\subsection{}The construction in
section \ref{assi} relies on the work of Fried, Haran, Jarden and
Sacerdote  on Galois stratifications
and 
quantifier elimination for pseudo-finite fields
\cite{F-S}, \cite{F-H-J}, \cite{F-J}.
More precisely,
quantifier elimination for pseudo-finite fields allows us to replace our
original  formula
by a quantifier free Galois formula, {\it i.e.}
a quantifier free formula in the language of rings and constants in
$k$,
extended by Artin symbols ${\rm Ar}_{W/U} (x)$ corresponding to
the decomposition subgroup of ${\rm Gal} (W / U)$ at $x$ modulo
conjugation,
with $x$ a variable running over the rational points of $U$. Here $W$ is
an {\'e}tale Galois cover of a normal algebraic variety $U$.
Then, by using Artin's
theorem on induced representation - a tool already present in the work of
Fried and Jarden 
\cite{F-J} -  and also
recent results on motives
\cite{G-S}, \cite{G-N}, \cite{B-N}, we are
able to associate a canonical virtual motive to
every quantifier free Galois formula.

\subsection{}In a previous paper \cite{Arcs}, we developped 
a general theory of integration - called motivic integration -
on the space of arcs of an algebraic variety $X$.
More precisely, let $X$ be an algebraic variety over a field $k$ of
characteristic zero. The space of arcs $\cL (X)$
is a $k$-scheme whose
$K$-rational points $\cL (X) (K)$ are exactly 
formal arcs $ {\rm Spec} \, K [[t]] \rightarrow X$,
for  every field $K$ containing $k$.
We also introduced in \cite{Arcs} an important boolean subalgebra
$\BB_{X}$ of the algebra of subsets of $\cL (X)$, the 
subalgebra of semi-algebraic subsets of $\cL (X)$. Semi-algebraic
subsets of $\cL (X)$ are definable in the language of valued fields.
We defined  a measure $\mu$ on $\BB_{X}$
with values - following ideas of Kontsevich - into the ring
$\widehat K_{0} ({\rm Sch}_{k})$ which is a certain completion
of the localization 
of the ring $K_{0} ({\rm Sch}_{k})$ with respect to the class of
the affine line.
In fact, there exists a general notion of measurable subset of 
$\cL (X)$ for $\mu$, developped in the appendix of \cite{MK},
and 
$\mu$ may be extended in a natural way to these measurable subsets.
See also 
\cite{Motivic} and 
\cite{TS} for variants with values into Grothendieck groups of motives.
This theory, which should maybe be called 
\textit{geometric} motivic integration, is not well suited for the aims of the 
present paper, the main reason being that at the residue
field
level, images by morphisms are considered
geometrically, in the sense of algebraic geometry,
and not at the level of looking to rational points,
which corresponds
to leaving the world of algebraic geometry  (polynomial equations or
inequations) for the world of first order logic of fields
(formulas with polynomials and quantifiers).
We are thus led to develop a different kind of 
motivic integration
theory, \textit{arithmetic} motivic integration,
which takes rationality properties in account.
This theory assigns to every definable subset of $\cL (X)$
an element 
in the ring
$\widehat K_{0}^{v} ({\rm Mot}_{k, \bar \QQ}) \otimes {\QQ}$,
where $\widehat K_{0}^{v} ({\rm Mot}_{k, \bar \QQ})$
is a certain completion of the localization of
$K_{0}^{v} ({\rm Mot}_{k, \bar \QQ})$ with respect to the Lefschetz motive.
Its very definition relies on  our construction of assigning a virtual
motive $\chi_{c} (\varphi)$ to a  formula
$\varphi$ in the first order language
of rings
with coefficients in $k$, so it involves in an essential
way arithmetical tools such as
decomposition subgroups, Chebotarev's Theorem and pseudo-finite fields.

\subsection{}Let us now review briefly
the content of the paper.
After a first section devoted to preliminaries, we develop in section
\ref{gs} what we need
on Galois stratifications. In section \ref{assi} we give our basic construction
of
assigning a virtual
motive to a  formula
in the first order language of rings.
First order formulas define subsets of the affine space. Since we need
to perform some geometric constructions like blowing ups, we 
introduce a new class of
objects, which are more geometric in nature
than formulas, which we call ``definable subassignements''
since they are not
in general functors.
Section \ref{dsr} is devoted to definable subassignements for rings and
section \ref{dspr} to definable subassignements for power series rings.
In section \ref{amids} we develop the basic theory of arithmetic motivic
integration.
Section \ref{rs} is devoted to general rationality results, in the spirit
of \cite{rat} and
\cite{Arcs}. We are then able in section \ref{interpol} to prove that 
arithmetic motivic integration specializes to $p$-adic
integration. In section \ref{versus} we define
the arithmetical Poincar{\'e} series
$P_{{\rm ar}} (T)$
for varieties $X$ over a field $k$ of characteristic zero
and we show it specializes to the $p$-adic Poincar{\'e} series when $k$ is
a number field.
The Poincar{\'e} series
$P_{{\rm ar}} (T)$
seems to contain much more interesting geometric information about the
variety $X$
than its geometric counterpart 
$P_{{\rm geom}} (T)$ introduced in \cite{Arcs}.
As an example, we compute both series for branches of plane 
curves in
section
\ref{ex}. In that
case the poles of $P_{{\rm ar}} (T)$ determine completely the Puiseux
pairs of the branch, while 
$P_{{\rm geom}} (T)$ allows one only to recover the multiplicity at the origin.

\setcounter{tocdepth}{1}

\tableofcontents

\section{Preliminaries}\label{prel}

\subsection{}For any ring $R$, by the first order language
of rings
with coefficients in $R$, we shall mean
the first order language (in the sense of logic),
with symbols to denote $+$, $-$, $\times$, $0$, $1$,
and for each element of $R$ a symbol to denote that element.
As for any first order language, formulas are built up from the above
symbols
together with variables, the logical connectives
$\wedge$ (and), $\vee$ (or), $\neg$ (not), the quantifiers $\exists$,
$\forall$ and the equality symbol $=$.

Let $S$ be a scheme. By a variety over $S$ or $S$-variety we shall
mean
a separated and reduced scheme
of finite type over $S$. When $S = {\rm Spec} \, R$ is affine we shall
also say variety over $R$ or $R$-variety.

\subsection{}\label{1.2}Let $k$ be a field. We shall
denote by $K_{0} ({\rm Sch}_k)$ the Grothendieck ring
of algebraic varieties over $k$.
It is the ring generated by symbols $[S]$, for $S$ an algebraic variety over
$k$, 
with the relations
$[S] = [S']$ if $S$ is isomorphic to $S'$,
$[S] = [S \setminus S'] + [S']$ if $S'$ is closed in $S$
and
$[S \times S'] = [S] \, [S']$.
Note that, for $S$ an algebraic variety
over
$k$,
the mapping $S' \mapsto [S']$
from the set of closed subvarieties of  $S$ extends uniquely to a mapping
$W \mapsto [W]$
from the set of constructible subsets of $S$ to $K_{0} ({\rm Sch}_k)$,
satisfying
$[W \cup W'] = [W] + [W'] - [W \cap W']$.
We set $\LL := [\AA^1_k]$ and
$K_{0} ({\rm Sch}_k)_{\rm loc} := K_{0} ({\rm Sch}_k) [\LL^{-1}]$.

Let $S$ be an algebraic variety over $k$.
We write $\dim S \leq n$
if all irreducible components of
$S$ have dimension $\leq n$.
Similarly, for $M$ in $K_{0} ({\rm Sch}_k)$,
we write $\dim M \leq n$ if
$M$ may be  expressed as a linear combination
of algebraic varieties with $\dim \leq n$.
For $m$ in $\ZZ$, we denote
by  $F^m K_{0} ({\rm Sch}_k)_{\rm loc}$
the subgroup of $K_{0} ({\rm Sch}_k)_{\rm loc}$ generated by
elements of the form
$[S] \, \LL^{- i}$ with 
$i - \dim S \geq m$. This
defines a decreasing filtration $F^m$
on $K_{0} ({\rm Sch}_k)_{\rm loc}$.
We denote by $\widehat K_{0} ({\rm Sch}_k)$ the completion of 
$K_{0} ({\rm Sch}_k)_{\rm loc}$ with respect to that filtration.
We do not know whether the natural morphism
$K_{0} ({\rm Sch}_k)_{\rm loc} \rightarrow
\widehat K_{0} ({\rm
Sch}_k)$ is injective or not.
We denote by $\overline K_{0} ({\rm Sch}_k)_{\rm loc}$ the image of 
$K_{0} ({\rm Sch}_k)_{\rm loc}$ in $\widehat K_{0} ({\rm Sch}_k)$.

\subsection{}\label{Chow}We denote by ${\rm Mot}_{k}$ the category of Chow motives
over $k$, with coefficients in $\QQ$. It is a pseudo-abelian
category and we denote by 
$K_{0} ({\rm Mot}_{k})$ its Grothendieck group.
It may also be defined as the abelian group associated to the monoid of
isomorphism
classes of objects in ${\rm Mot}_{k}$ with respect to $\oplus$.
Let us recall (see \cite{Sc} for more details),
that objects in ${\rm Mot}_{k}$ are just triples
$(S, p, n)$ with $S$ proper and smooth over $k$, $p$ an idempotent
correspondence with coefficients in $\QQ$
on $S$ and $n$ in $\ZZ$.
The tensor
product on ${\rm Mot}_{k}$ induces a product on $K_{0} ({\rm Mot}_{k})$
which provides $K_{0} ({\rm Mot}_{k})$ with 
a natural ring structure.

Assume now that $k$ is of characteristic zero.
By a result of Gillet and Soul{\'e} \cite{G-S}
and Guill{\'e}n and Navarro Aznar \cite{G-N} there exists a unique
morphism of rings
$$
\chi_{c} : K_{0} ({\rm Sch}_k) \longrightarrow K_{0} ({\rm Mot}_{k})
$$
such that $\chi_{c} ([S]) = [h (S)]$
for $S$ projective and smooth, where $h (S)$ denotes the Chow motive
associated to $S$, {\it i.e.} $h (S) = (S, {\rm id}, 0)$.
Let us still denote by $\LL$ the image of $\LL$ by $\chi_{c}$.
Since $\LL = [({\rm Spec} \, k, {\rm id}, -1)]$, it is invertible in
$K_{0} ({\rm Mot}_{k})$, hence $\chi_{c}$ can be extended uniquely to
a ring morphism
$$
\chi_{c} : K_{0} ({\rm Sch}_k)_{\rm loc} \longrightarrow K_{0} ({\rm Mot}_{k}).
$$

Let us denote by ${\rm Mot}_{k, \bar \QQ}$ the category of Chow motives
over $k$, with coefficients in $\bar \QQ$
and by
$K_{0} ({\rm Mot}_{k, \bar \QQ})$ its Grothendieck group.
Objects in ${\rm Mot}_{k, \bar \QQ}$ are triples
$(S, p, n)$ with $S$ proper and smooth over $k$, $p$ an idempotent
correspondence with coefficients in $\bar \QQ$
on $X$ and $n$ in $\ZZ$,
see, {\it e.g.,}
\cite{Motivic} for more details.
We denote by
$K_{0}^{v} ({\rm Mot}_{k, \bar \QQ})$,
resp.  $K_{0}^{v} ({\rm Mot}_{k, \bar \QQ})_{\rm loc}$,
the image of
$K_{0} ({\rm Sch}_k)$, resp. $K_{0} ({\rm Sch}_k)_{\rm loc}$, in 
$K_{0} ({\rm Mot}_{k, \bar \QQ})$ by the morphism,
which we will still denote by $\chi_{c}$, obtained by
composition of $\chi_{c}$ with
the natural morphism 
$K_{0} ({\rm Mot}_{k}) \rightarrow K_{0} ({\rm Mot}_{k, \bar \QQ})$.
We denote by $F^m$ the decreasing filtration
on $K_{0}^{v} ({\rm Mot}_{k, \bar \QQ})_{\rm loc}$ which is the image of the
filtration $F^{m}$ on $K_{0} ({\rm Sch}_k)_{\rm loc}$,
and by 
$\widehat K_{0}^{v} ({\rm Mot}_{k, \bar \QQ})$ the
completion of
$K_{0}^{v} ({\rm Mot}_{k, \bar \QQ})_{\rm loc}$ with respect to the filtration
$F^{m}$.
We also define  
$\overline K_{0}^{v} ({\rm Mot}_{k, \bar \QQ})_{\rm loc}$ as the image of
$K_{0}^{v} ({\rm Mot}_{k, \bar \QQ})_{\rm loc}$
in $\widehat K_{0}^{v} ({\rm Mot}_{k, \bar \QQ})$.
We set 
\begin{gather*}K_{0}^{v} ({\rm Mot}_{k, \bar \QQ})_{\QQ}~:=
K_{0}^{v} ({\rm Mot}_{k, \bar \QQ}) \otimes \QQ,
 K_{0}^{v} ({\rm Mot}_{k, \bar \QQ})_{{\rm loc}, \QQ}~:=
 K_{0}^{v} ({\rm Mot}_{k, \bar \QQ})_{\rm loc} \otimes \QQ,
\\
\overline K_{0}^{v} ({\rm Mot}_{k, \bar \QQ})_{{\rm loc}, \QQ}~:=
\overline K_{0}^{v} ({\rm Mot}_{k, \bar \QQ})_{\rm loc} \otimes \QQ
\,
{\rm and}
\,
\widehat K_{0}^{v} ({\rm Mot}_{k, \bar \QQ})_{\QQ}~:=
\widehat K_{0}^{v} ({\rm Mot}_{k, \bar \QQ}) \otimes \QQ.
\end{gather*}
We denote by $F^{m} \widehat K_{0}^{v} ({\rm Mot}_{k, \bar \QQ})_{\QQ}$
the filtration on $\widehat K_{0}^{v} ({\rm Mot}_{k, \bar
\QQ})_{\QQ}$ naturally
induced by
$F^{m}$ .

If $k'$ is a field containing $k$, and $M$ belongs to
${\rm Mot}_k$ we shall denote by $M \otimes k'$ the same object but
considered as an element of ${\rm Mot}_{k'}$. Similarly, one denotes by
$M \otimes k'$ the image of an element $M$ in
any of the preceding rings relative to $k$
in the similar ring relative to $k'$.

\section{Galois stratifications}\label{gs}

\subsection{Galois stratifications}\label{GS}
Let $A$ be an integral and normal scheme.
A morphism of schemes
$h : C \rightarrow A$ is a
{\it Galois cover} if $C$ is integral,
$h$ is {\'e}tale (hence $C$ is normal)
and there is a
finite group $G = G (C  / A)$, the Galois group, acting on $C$ 
such that $A$  is isomorphic to the quotient $C / G$ 
in such a way that
$h$ is the composition of the quotient morphism with the 
isomorphism. Isomorphisms of Galois covers are defined
in the usual way.
If $A'$ is a locally closed integral and normal subscheme
of $A$, let $C'$ denote any connected component 
of $C \cap h^{-1} (A')$. One defines the restriction of the Galois cover
$h : C \rightarrow A$ to $A'$ as the Galois cover 
$h' : C' \rightarrow A'$, with $h'$ the restriction of $h$ to 
$C'$. The Galois group
$G (C' / A')$ is the decomposition
subgroup of $G (C / A)$ at the generic point of $C'$.
The choice of another connected component would give an isomorphic
Galois cover.
We say that the Galois cover $h : C \rightarrow A$
is \textit{colored}
if $G (C  / A)$ is equipped with
a family ${\rm Con}$  of subgroups of
$G (C  / A)$ 
which is stable by conjugation under elements of
$G (C  / A)$.
The restriction of a colored  Galois cover
$h : C \rightarrow A$ to $A'$,
a locally closed integral and normal subscheme
of $A$, is defined by the family ${\rm Con}'$  of subgroups of
$G (C'  / A')$ which belong to ${\rm Con}$.

Let $S$ be an integral normal scheme and let $X_{S}$ be a variety
over $S$.
A {\it normal stratification}
of $X_{S}$, 
$$
\cA = \langle X_{S}, C_{i} / A_{i} \, \vert \, i \in I \rangle , 
$$
is a partition of $X_{S}$
into a finite
set of integral and normal locally closed $S$-subschemes
$A_{i}$, $i \in I$, each equipped with a 
Galois cover $h_{i} : C_{i} \rightarrow A_{i}$.
A normal stratification $$\cA' = \langle X_{S},
C'_{i} / A'_{i} \, \vert \, i' \in I' \rangle
$$ will be said to be {\it finer} than
$\cA$ if, for each $i$ in $I'$,
$A'_{i}$ is contained in some
$A_{j}$ and $C'_{i} \rightarrow A'_{i}$ is isomorphic
to
the restriction of $C_{j} \rightarrow A_{j}$
to  $A'_{i}$ as a Galois cover. We will say that
normal stratifications $\cA$ and $\cA'$ are isomorphic if
there exists a normal stratification $\cA''$ which is finer than both
$\cA$ and $\cA'$.

A normal stratification $\cA$
may be augmented to a {\it Galois stratification} 
$$
\cA 
= \langle X_{S}, C_{i} / A_{i}, {\rm Con} (A_{i})
\, \vert \, i \in I \rangle, 
$$
if for each $i \in I$, ${\rm Con} (A_{i})$
is a family of subgroups of
$G (C_{i}  / A_{i})$ 
which is stable by conjugation under elements in $G (C_{i}  / A_{i})$,
{\it i.e.} $(C_{i} \rightarrow A_{i}, {\rm Con} (A_{i}))$
is a colored Galois cover.
One defines similarly as before finer  and isomorphic
Galois stratifications.
In general it will be harmless to identify 
isomorphic Galois stratifications.

To any $S$-constructible subset $W$ of
$X_{S}$ one may associate a well defined (up to
isomorphism) Galois stratification by taking any normal stratification
with all strata contained either in $W$ or in its complement,
by taking the identity morphism as Galois cover on each stratum, 
and taking for ${\rm Con} (A_{i})$ the family consisisting of
the group with one element 
if
$A_{i} \subset W$ and the empty family otherwise.

We define the support ${\rm Supp} \,  (\cA) $
of a Galois stratification
$$\cA 
= \langle X_{S}, C_{i} / A_{i}, {\rm Con} (A_{i})
\, \vert \,  i \in I \rangle, 
$$
as the union of the sets $A_{i}$ with ${\rm Con} (A_{i})$ non empty. We 
define
the dimension ${\rm dim} \, (\cA)$ as the maximum of the dimensions 
of the sets $A_{i}$ with ${\rm Con} (A_{i})$ non empty.

\subsection{Galois formulas}\label{Galois formulas}Let $U = {\rm Spec} \, R$
be an affine scheme, which we assume to be integral and normal.
For any variety $X_{U}$ over $U$
and any closed point $x$
of $U$, we denote by $\FF_{x}$ the residual field of $x$ on $U$
and by
$X_{x}$ the fiber of $X_{U}$ at $x$.
More generally, for any field $M$ containing $\FF_{x}$, we shall denote by
$X_{x, M}$ the fiber product of $X_{U}$ and ${\rm Spec} \, M$ over $x$.

Let $X_{U}$ be a variety over $U$.
Let $\cA 
= \langle X_{U}, C_{i} / A_{i}, {\rm Con} (A_{i}) \, \vert \,
i \in I \rangle$
be a Galois stratification of $X_{U}$ and let
$x$ be a closed point of $U$.
Let $M$ be a field containing $\FF_{x}$ and 
let $\aa$
be  an $M$-rational point of $X_{x}$ belonging
to $A_{i, x}$. We denote by ${\rm Ar} (C_{i} / A_{i}, x ,\aa)$
the conjugacy class of subgroups of $G (C_{i}  / A_{i})$
consisting of  the  decomposition subgroups at $\aa$. We shall
write
$$
{\rm Ar} (\aa) \subset {\rm Con} (\cA)
$$
for 
$${\rm Ar} (C_{i} / A_{i}, x ,\aa) \subset {\rm Con} (A_{i}).
$$
For $x$ a  closed point of $U$ and
$M$ a field containing $\FF_{x}$, we consider the subset
$$
Z (\cA, x, M) :=
\Bigl\{\aa \in X_{U} (M)
\Bigm \vert {\rm Ar} (\aa) \subset
{\rm Con} (\cA)
\Bigr\}
$$
of
$X_{U}(M)$.

Let
$\cA 
= \langle \AA^{m + n}_{U}, C_{i} / A_{i}, {\rm Con} (A_{i})
\, \vert \, i \in I \rangle$
be a Galois stratification of $\AA^{m + n}_{U}$.
Let $Q_{1}, \ldots, Q_{m}$ be quantifiers.
We denote by $\vartheta$  or by $\vartheta (\YY)$ the formal expression
$$
(Q_{1} X_{1}) \ldots (Q_{m} X_{m}) \, [{\rm Ar} (\XX, \YY) \subset
{\rm Con} (\cA)]
$$
with $\XX = (X_{1}, \ldots, X_{m})$ and
$\YY = (Y_{1}, \ldots, Y_{n})$. We call $\vartheta (\YY)$
a {\it Galois formula} over $R$ in the free variables $\YY$.

Now to a Galois formula $\vartheta$, to a 
closed point $x$ of $U$ and to a field $M$ containing $\FF_{x}$,
one associates the subset
\begin{multline*}
Z (\vartheta, x, M) :=\\
\Bigl\{\bb = (b_{1}, \ldots, b_{n}) \in M^n
\Bigm \vert
(Q_{1} a_{1}) \ldots (Q_{m} a_{m}) \, [{\rm Ar} (\aa, \bb) \subset
{\rm Con} (\cA)]
\Bigr\}
\end{multline*}
of
$M^n$, where the quantifiers $Q_{1} a_{1}$, \dots, $Q_{m} a_{m}$ run
over $M$.

Let $\varphi (Y_{1}, \ldots, Y_{n})$ be a formula 
in the first order
language of rings with coefficients in the ring $R$ and free variables
$Y_{1}, \ldots, Y_{n}$.
For every 
closed point $x$ in $U$ and every field $M$ containing
$\FF_{x}$ we denote by 
$Z (\varphi, x, M)$ the subset of 
$M^n$ defined by the (image over $\FF_{x}$ of the) formula
$\varphi$.
Assume now $\varphi$ is
in prenex normal form,
{\it i.e.} a formula of the form
\begin{equation}\label{prenex}
(Q_{1} X_{1}) \ldots (Q_{m} X_{m}) \, \Bigl[
\bigvee_{i = 1}^{k} \bigwedge_{j = 1}^l f_{i, j}
(\XX, \YY) = 0 \wedge g_{i, j}
(\XX, \YY) \not= 0
\Bigr],
\end{equation}
with $f_{i, j}$ and $g_{i, j}$ in $R [\XX, \YY]$.
The formula between brackets defines
an $U$-construct\-ible subset $A$ of $\AA_{U}^{m + n}$ to which one 
associates a Galois stratification as above. In this way
one obtains a Galois formula $\vartheta$ over $R$ such that
$Z (\vartheta, x, M) = Z (\varphi, x, M)$ for every closed point $x$ in
$U$ and every field $M$
containing $\FF_{x}$.

\subsection{Quantifier elimination for Galois formulas}

We recall that a pseudo-finite field
$F$ is a perfect
infinite field which has exactly one extension of each degree
and such that every absolutely irreducible variety over $F$
has a rational point
over
$F$. J. Ax proved \cite{Ax} that an infinite field $F$ is pseudo-finite
if and only if every sentence\footnote{A sentence is a formula without free
variables} in the first order language
of rings which is true in all finite fields is also true in $F$.
We recall also that the property of being a pseudo-finite field
is stable by ultraproducts\footnote{see, {\it e.g.}, \cite{F-J} for the
definition of ultraproducts}.

There exists 
several versions of quantifier elimination for Galois formulas
\cite{F-S}, \cite{F-H-J}, \cite{F-J}.
The following one seems to be best suited for the present work.

\begin{theorem}\label{pfelimination}Let $k$ be a field.
Let $\cA$
be a Galois stratification of $\AA^{m + n}_{k}$ and let  $\vartheta$  
be a
Galois formula
$$
(Q_{1} X_{1}) \ldots (Q_{m} X_{m}) \, [{\rm Ar} (\XX, \YY) \subset
{\rm Con} (\cA)]
$$
with respect to $\cA$.
There exists a  Galois stratification
$\cB$ of 
$\AA^{n}_{k}$ such that,
for every pseudo-finite field $F$  containing $k$,
$$Z (\vartheta, {\rm Spec} \, k, F)
=
Z (\cB, {\rm Spec} \, k, F).
$$
\end{theorem}

\begin{proof}Since pseudo-finite fields are Frobenius fields in the
terminology of \cite{F-J}, the result is a special case
of Proposition 25.9 of \cite{F-J}.
\end{proof}

\begin{cor}\label{pf2.5}Let $\varphi (Y_{1}, \ldots, Y_{n})$ be a formula 
in the first order
language of rings with coefficients in a field $k$ and free variables
$Y_{1}, \ldots, Y_{n}$.
There exists a Galois stratification
$\cB$ of 
$\AA^{n}_{k}$ such that, 
for every pseudo-finite  field $F$ containing $k$,
$$Z (\varphi, {\rm Spec} \, k, F)
=
Z (\cB, {\rm Spec} \, k, F).
$$
\end{cor}

\begin{proof}Take $\varphi'$ a formula in prenex normal form,
which is logically equivalent to $\varphi$. Since,
for every pseudo-finite  field $F$ containing $k$,
$Z (\varphi, \, k, F) = Z (\varphi', \, k, F)$, the result follows from 
Theorem \ref{pfelimination} and the observation made
at the end of \ref{Galois formulas}.
\end{proof}

\begin{remark}It follows from \cite{F-J} that there is an effective 
algorithm to determine a
Galois stratification
$\cB$  in Theorem \ref{pfelimination} and Corollary \ref{pf2.5}.
\end{remark}

\subsection{}Let $k$ be a field
and let $\varphi$ and $\varphi'$
be two formulas in the first order language
of rings with coefficients in $k$
and free variables $X_{1}, \cdots,
X_{n}$. Let $F$ be a field containing $k$.
Recall that one says that
$\varphi$ and $\varphi'$ are equivalent in $F$ if they define the same
subsets in $F^{n}$.
We define the equivalence relation
$\approx$ by $\varphi \approx \varphi'$ if
$\varphi$ is equivalent to $\varphi'$ in every pseudo-finite field $F$
containing $k$.
There is also a weaker equivalence relation $\equiv$ defined as follows.
Let $\varphi$ and $\varphi'$ two formulas in the first order language
of fields with coefficients in $k$
and free variables $X_{1}, \cdots,
X_{n}$ and $X_{1}, \cdots,
X_{n'}$, respectively. We write $\varphi \equiv \varphi'$ if there exists
a formula $\psi$ in the first order language
of rings with coefficients in $k$
and free variables $X_{1}, \cdots,
X_{n + n'}$, such that, for every 
pseudo-finite field $F$
containing $k$, the subset $Z (\psi, \spec k, F)$ of $F^{n + n'}$ is the
graph of a bijection between 
$Z (\varphi, \spec k, F)$ and $Z (\varphi', \spec k, F)$.
In a few cases we shall write $\approx_{k}$ and $\equiv_{k}$ to keep
track of
the field $k$.

To get a more concrete interpretation of the equivalence relations
$\approx$ and $\equiv$,
we shall assume now that $k$ is a field of characteristic zero
which is the field of fractions of a normal domain $R$
of finite type over $\ZZ$.
We set $U = {\rm Spec} \, R$, and,
for any non zero element $f$ of $R$, we set
$U_{f} = {\rm Spec} \, R [f^{-1}]$.
Let $\cA_{U}$ be a normal (resp. Galois) stratification of $\AA^n_{U}$.
By base change  over $\Spec k$,
one
obtains a normal (resp. Galois) stratification
of $\AA^n_{k}$, which we will denote by $\cA_{U} \otimes k$.
Conversely, any normal or Galois stratification of
$\AA^n_{k}$ may be obtained in this way, at the cost of replacing $U$ by
some localization $U_{f}$.

\begin{prop}\label{star}\begin{enumerate}
\item[(1)]Let $\varphi$ and $\varphi'$
be two formulas in the first order language
of rings with coefficients in $k$
and free variables $X_{1}, \cdots,
X_{n}$. Viewing $\varphi$ and $\varphi'$ as  formulas
in the first order language
of rings with coefficients in $R_{f}$, for a suitable $f$,
we have
$\varphi \approx \varphi'$ if and only if,
for a suitable non zero multiple $f'$ of $f$, 
$Z (\varphi, x, \FF_{x}) = Z (\varphi', x,
\FF_{x})$,
for every closed point $x$ of
$U_{f'}$.
\item[(2)]Let $\varphi$ and $\varphi'$ two formulas in the first order language
of rings with coefficients in $k$
and free variables $X_{1}, \cdots,
X_{n}$ and $X_{1}, \cdots,
X_{n'}$, respectively. 
Viewing $\varphi$ and $\varphi'$ as  formulas
in the first order language
of rings with coefficients in $R_{f}$, for a suitable $f$,
we have
$\varphi \equiv \varphi'$ if and only if, for a suitable non zero
multiple $f'$ of $f$, 
there exists
a formula $\psi$ in the first order language
of rings with coefficients in $R_{f'}$
and free variables $X_{1}, \cdots,
X_{n + n'}$, such that,
for every closed point $x$ of
$U_{f'}$,
the subset $Z (\psi, x, \FF_{x})$ of $\FF_{x}^{n + n'}$ is the
graph of a bijection between 
$Z (\varphi, x, \FF_{x})$ and $Z (\varphi', x, \FF_{x})$.
\end{enumerate}
\end{prop}

\begin{proof}This  follows from Proposition \ref{2star}
applied to the sentences
$$\forall x (\varphi(x) \longleftrightarrow \varphi' (x))$$
and
$$
[\forall x (\varphi (x) \rightarrow \exists! x' : (\varphi' (x')
\wedge \psi (x, x')))]
\wedge
[\forall x' (\varphi' (x') \rightarrow \exists! x : (\varphi (x)
\wedge \psi (x, x')))],
$$
in case (1) and (2), respectively.
\end{proof}

\begin{prop}\label{2star}Let $k$ be a field of characteristic zero
which is the field of fractions of a normal domain $R$
of finite type over $\ZZ$. Let $\vartheta$ be a Galois formula
over $R$ (as in \ref{Galois formulas}) with no free variables.
Then $\vartheta$ is true in every pseudo-finite field $F$ containing
$k$ if and only if there exists $f$ in $R \setminus \{0\}$,
such that, for every closed point $x$ in ${\rm Spec} \, R_{f}$,
the formula $\vartheta$ is true in $\FF_{x}$,
{\it i.e.} $Z (\vartheta, x, \FF_{x}) = \AA^{0}_{\FF_{x}}(\FF_{x})$.
\end{prop}

\begin{proof}We first assume that $\vartheta$ is true in every
pseudo-finite field $F$ containing
$k$. If, for every $f$ in $R \setminus \{0\}$, there would exist
a closed point $x$ in ${\rm Spec} \, R_{f}$ such that $\vartheta$ is
false in $\FF_{x}$, a suitable ultraproduct of
the fields $\FF_{x}$ would yield a pseudo-finite field
containing $k$ in which $\vartheta$ is false,
since the ultraproduct construction commutes in the
present case with the Artin symbol.

Conversely, suppose that there exists $f$ in $R \setminus \{0\}$
such that, for every closed point $x$ in ${\rm Spec} \, R_{f}$,
the formula $\vartheta$ is true in $\FF_{x}$.
By the quantifier elimination Theorem \ref{pfelimination},
there exists a Galois formula $\cB$, over a suitable $R_{f}$,
with no free variables and no quantifiers, such that $\vartheta
\leftrightarrow
\cB$ holds in every pseudo-finite field $F$ containing $k$.
It follows from the first part of the proof that, maybe only after
replacing $f$ by a multiple,
$\vartheta
\leftrightarrow
\cB$ holds
also in $\FF_{x}$ for every closed point $x$ in ${\rm Spec} \, R_{f}$.

Thus we may suppose that $\vartheta$ has no quantifiers. Let $\cA$ be
the Galois stratification over $R$ belonging to $\vartheta$. Because
$\vartheta$
has no free variables, $\cA \otimes k$ consists of only one cover
${\rm Spec} \, L \rightarrow {\rm Spec} \, k$, with $L$ a field which is
Galois over $k$.

Assume now that there exists a pseudo-finite field $F$ containing $k$ in
which $\vartheta$ is false. Let $\sigma$ be a topological generator
of the absolute Galois group of $F$, and denote by  $\tau$
the restriction of $\sigma$ to $L$. Then 
$$
{\rm Ar} \, ({\rm Spec} \, L / {\rm Spec} \, k, 0, F)
=
\Bigl\{ \alpha <\tau> \alpha^{-1} \Bigm\vert \alpha \in {\rm Gal} \, (L
/ k)
\Bigr\},
$$
where $<\tau>$ denotes the subgroup of
${\rm Gal} \, (L
/ k)$ generated by $\tau$.
Since we suppose that $\vartheta$ is false in $F$, we have
$<\tau> \not\subset {\rm Con} (\cA) $. By Chebotarev's Theorem
for ${\rm Spec} \, R$ (see, {\it e.g.}, \cite{Cheb}), there exists a closed point $x$ of
${\rm Spec} \, R_{f}$ such that
$$
{\rm Ar} \, ({\rm Spec} \, L / {\rm Spec} \, k, x, 0)
=
\Bigl\{ \alpha <\tau> \alpha^{-1} \Bigm\vert \alpha \in {\rm Gal} \, (L
/ k)
\Bigr\}.
$$
Thus 
$
{\rm Ar} \, ({\rm Spec} \, L / {\rm Spec} \, k, x, 0)
\not \subset {\rm Con} (\cA) $, whence $\vartheta$ is false in $\FF_{x}$,
which contradicts our assumption that $\vartheta$ is true in $\FF_{x}$.
\end{proof}

\begin{lem}\label{rest}Let $\varphi$ and $\varphi'$
be two formulas in the first order language
of rings with coefficients in a field $k$.
If $\varphi \approx_{k} \varphi'$ (resp. $\varphi
\equiv_{k} \varphi'$), then there exists a subfield
$k_{0}$ of $k$, of finite type over its prime field
and such that the coefficients in $\varphi$ and $\varphi'$ belong to
$k_{0}$, such that $\varphi \approx_{k_{0}} \varphi'$
(resp. $\varphi
\equiv_{k_{0}} \varphi'$).
\end{lem}

\begin{proof}First consider the case where $\varphi \approx_{k} \varphi'$. 
Let $\{k_{i} \, \vert \, i \in \Sigma \}$
be the set of all subfields of $k$
which are of finite type over the prime field of $k$ and which
contain the
constants in $\varphi$ and $\varphi'$, and assume there
exists for each $i$ in $\Sigma$ a pseudo-finite field $F_{i}$ containing
$k_{i}$ in which $\varphi$ is not equivalent to $\varphi'$. Choose an
ultrafilter $D$ on $\Sigma$ containing $S_{i} := \{j \, \vert \, k_{i}
\subset k_{j}\}$, for every $i $ in $\Sigma$. Let $F$ be the
ultraproduct of all the fields
$F_{i}$ with respect to $D$. Clearly $F$ is a
pseudo-finite field in which  
$\varphi$ is not equivalent to $\varphi'$. Moreover $k$ is imbedded in
$F$ by the map
$a \mapsto (a_{i})_{i \in I} \mod D$, where $a_{i} = a$ if $a \in k_{i}$
and $a_{i} = 0$ if $a \not\in k_{i}$. This contradicts $\varphi \approx
\varphi'$. The proof for $\equiv$ is just the same,
considering now the set $\{k_{i} \, \vert \, i \in \Sigma \}$
of all subfields of $k$
which are of finite type over the prime field of $k$ and contain the
constants in $\varphi$, $\varphi'$ and $\psi$, where $\psi$ is a formula
satisfying the conditions in the definition of $\equiv$ and replacing
everywhere
``$\varphi$ is not equivalent to $\varphi'$ in the pseudo-finite field $M$''
by
``$Z (\psi, \spec k, M)$ is not the
graph of a bijection between 
$Z (\varphi, \spec k, M)$ and $Z (\varphi', \spec k, M)$''.\end{proof}

The following quantifier elimination statement follows directly from
Theorem \ref{pfelimination} and Proposition \ref{2star}.

\begin{prop}\label{elimination}
Let $U = {\rm Spec} \, R$
be an affine scheme of finite 
type over $\ZZ$, integral, normal and of characteristic zero.
Let $\cA$
be a Galois stratification of $\AA^{m + n}_{U}$ and let  $\vartheta$  
be a
Galois formula
$$
(Q_{1} X_{1}) \ldots (Q_{m} X_{m}) \, [{\rm Ar} (\XX, \YY) \subset
{\rm Con} (\cA)]
$$ with $Q_{1}, \ldots, Q_{m}$ quantifiers.
There exists a non zero element $f$ of $R$ and a Galois stratification
$\cB$ of 
$\AA^{n}_{U_{f}}$ such that,
for every closed point $x$ of
$U_{f}$,
$$Z (\vartheta, x, \FF_{x})
=
Z (\cB, x, \FF_{x}).\qed
$$
\end{prop}

\begin{cor}\label{2.5}Let $\varphi (Y_{1}, \ldots, Y_{n})$ be a formula 
in the first order
language of rings with coefficients in the ring $R$ and free variables
$Y_{1}, \ldots, Y_{n}$.
There exists a non zero element $f$ of $R$ and a Galois stratification
$\cB$ of 
$\AA^{n}_{U_{f}}$ such that,  for every closed point $x$ of
$U_{f}$,
$$Z (\varphi, x, \FF_{x})
=
Z (\cB, x, \FF_{x}).\qed
$$
\end{cor}

\section{Assigning virtual motives to formulas}\label{assi}

\subsection{Motives and group action}Let
$k$ be a field of characteristic zero.
Let $G$ be a finite group. 
Let $X$ be an algebraic variety over $k$ endowed with a $G$-action. 
We say $X$ is a $G$-variety if the $G$-orbit of every closed point in $X$
is contained in an affine open subscheme (this condition is
always satisfied when
$X$ is quasi-projective). One defines in the usual way isomorphisms 
and closed immersions of $G$-varieties and so one may define
a ring $K_{0} ({\rm Sch}_k, G)$, the Grothendieck ring
of $G$-varieties over $k$, similarly as in \ref{1.2}.

Let $\alpha$ be the character of a representation
$G \rightarrow {\rm GL} (V_{\alpha})$
defined over $\bar \QQ$. Denote by $n_{\alpha}$ the
dimension of $V_{\alpha}$ and consider
$$
p_{\alpha} := \frac{n_{\alpha}}{|G|} \sum_{g \in G} \alpha^{-1} (g) [g],
$$
the 
corresponding
idempotent in $\bar \QQ [G]$.

There is a natural ring morphism $\mu$
from 
$\bar \QQ [G]$ to the ring of correspondences on $X$ with coefficients
in $\bar \QQ$ sending a group element $g$
onto the graph of multiplication by $g$.

We will need the following equivariant version of
the result of Gillet and Soul{\'e} \cite{G-S}
and Guill{\'e}n and Navarro Aznar \cite{G-N}.
We denote by $R_{\bar \QQ} (G)$ 
the group of  characters of virtual representations of 
$G$
defined over $\bar \QQ$.

\begin{theorem}\label{twisted}Let $G$ be a finite group. For every virtual
character
$\alpha$ in $R_{\bar \QQ} (G)$, there exists a unique
morphism of rings
$$
\chi_{c} ( \_, \alpha) :
K_{0} ({\rm Sch}_k, G) \longrightarrow K_{0} ({\rm Mot}_{k, \bar \QQ})
$$
such that:
\begin{enumerate}\item[(1)]If  $X$  is
projective and smooth with $G$-action and $\alpha$ is the character of
an
irreducible representation defined over $\bar \QQ$, 
$ n_{\alpha} \chi_{c} ([X], \alpha)$ is equal to the class of the motive
$(X, \mu (p_{\alpha}), 0)$ in $K_{0} ({\rm Mot}_{k, \bar \QQ})$.
\item[(2)]For every $G$-variety $X$, 
$$
\chi_{c} (X) = \sum_{\alpha} n_{\alpha} \chi_{c} (X, \alpha)
$$
where $\alpha$ runs over all irreducible characters of $G$.
\item[(3)]For every $G$-variety $X$, the function
$\alpha \mapsto \chi_{c} (X, \alpha)$ is
a group morphism $R_{\bar \QQ} (G) \rightarrow 
K_{0} ({\rm Mot}_{k, \bar \QQ})$.
\end{enumerate}
\end{theorem}

\begin{proof}This is Theorem 6.1 of \cite{B-N}. When $G$ is abelian
it is Theorem 1.3.1 of
\cite{Motivic} where the hypothesis $G$ abelian
was not used seriously in the proof.
\end{proof}

We shall also use the following result.

\begin{prop}\label{BN}Let $G$ be a finite group, let $H$ be a subgroup of
$G$ and let $X$ be 
a  $G$-variety.
\begin{enumerate}\item[(1)]Assume $H$ is
a normal subgroup of $G$. Then, for every character 
$\alpha$ of $G / H$,
$$\chi_{c} (X / H, \alpha) = \chi_{c} (X, \alpha \circ \varrho),$$
where $\varrho$ is the projection $ G \rightarrow G / H$.
\item[(2)]Let $\alpha$ be a character of $H$. Then
$$
\chi_{c} (X, {\rm Ind}^{G}_{H} \alpha) = 
\chi_{c} (X, \alpha),
$$
viewing $X$ as an $H$-variety
in the second term of the equality.
\item[(3)]Assume $X$ is isomorphic as a $G$-variety to
$\bigoplus_{s \in G / H} s Y$, with $Y$ an $H$-variety.
Then, for every character 
$\alpha$ of $G$,
$$
\chi_{c} (X, \alpha) = \chi_{c} (Y, {\rm Res}_{H}^{G} \alpha).
$$
\end{enumerate}
\end{prop}

\begin{proof}Statements (1) amd (3) follow from Proposition 6.3 of 
\cite{B-N} 
applied to  the morphism $\varrho : G \rightarrow G/H$ and
$H \rightarrow G$, respectively.
Similarly, statement (2) follows from a dual form of 
Proposition 6.3 of 
\cite{B-N} which is as follows (notation being as in loc. cit.)~:
for $\psi : G \rightarrow G'$ a morphism of finite groups,
$\alpha$ a character of $G$ and $X$ a $G'$-variety,
we have $\chi_{c} ({\rm Res}_{\psi} X, \alpha) = \chi_{c} (X,
{\rm Ind}_{\psi} \alpha) $.
The proof is just similar to  the one of Proposition 6.3 of 
\cite{B-N} using the projection formula
Corollary 4.3 of \cite{B-N} and the fact 
that  the functor $h_{c}$ of loc. cit. commutes with ${\rm Res}_{\psi}$.
Another way to  prove (2) and (3) would be to remark that
when $X$ is projective and
smooth it is just
a consequence of Theorem \ref{twisted} (1) together with elementary
theory of representations of finite groups and then to deduce the result for
arbitrary $X$
by additivity of $\chi_{c}$.
\end{proof}

Let us denote by $C (G, \bar \QQ)$ the $\bar \QQ$-vector space
of  $\bar \QQ$-valued central functions on $G$. For every $W$
in $K_{0} ({\rm Sch}_k, G)$, one defines by linearity
a morphism of vector spaces
$\alpha \mapsto \chi_{c} (W,  \alpha)$
from $C (G, \bar \QQ)$
to $K_{0} ({\rm Mot}_{k, \bar \QQ}) \otimes \bar \QQ$,
expressing a central function as a linear
combination of irreducible characters.
Now let us denote by $C (G, \QQ)$ the $\QQ$-vector space of $\QQ$-central
functions $G \rightarrow  \QQ$, {\it i.e.} $\QQ$-linear combinations of
characters of 
$\QQ$-irreducible representations
of $G$ defined over $\QQ$.
We recall that a central function $\alpha : G \rightarrow  \QQ$
belongs to $C (G, \QQ)$ if and only if $\alpha (x) = \alpha (x')$
for each $x$, $x'$ in $G$ such  that $\langle x \rangle$ is conjugate
to $\langle x' \rangle$.
Here $\langle x \rangle$ denotes the subgroup generated by $x$.

\begin{prop}\label{qval}Let $G$ be a finite group. For every
central function $\alpha: G \rightarrow \QQ$ in $C (G, \QQ)$
and every $W$ in 
$K_{0} ({\rm Sch}_k, G)$, the virtual motive
$\chi_{c} (W,  \alpha)$ belongs to
the image of the morphism of
$\QQ$-vector spaces
$$K_{0} ({\rm Sch}_k) \otimes \QQ
\longrightarrow K_{0} ({\rm Mot}_{k, \bar \QQ}) \otimes \QQ$$
induced by $\chi_{c}$, {\it i.e.} is an element of
$K_{0}^{v} ({\rm Mot}_{k, \bar \QQ})_{\QQ}$.
\end{prop}

\begin{proof}By a classical result of Emil Artin, $\alpha$ is a
$\QQ$-linear combination of characters
of the form ${\rm Ind}^G_{H} 1_{H}$ with $H$ a cyclic subgroup
of $G$. It follows from Proposition \ref{BN} that, for every $G$-variety
$X$, $$\chi_{c} (X,  {\rm Ind}^G_{H} 1_{H})
= \chi_{c} (X, 1_{H}) =
\chi_{c} (X / H),$$
where in the middle term $X$ is viewed as an $H$-variety,
whence $\chi_{c} (X,  {\rm Ind}^G_{H} 1_{H})$ belongs to 
$K_{0}^{v} ({\rm Mot}_{k, \bar \QQ})$.
\end{proof}

\subsection{Assigning virtual motives to Galois stratifications}Let $k$ be a
field 
of characteristic zero and
let $A$ be a normal integral variety over $k$.
Let $h : C \rightarrow A$
be a Galois cover with Galois group $G$ and let
${\rm Con}$ be a family of subgroups of
$G$ 
which is stable by conjugation under elements in $G$.
One may associate to these
data an element
$\chi_{c} (C / A, {\rm Con})$ of
$K_{0}^{v} ({\rm Mot}_{k, \bar \QQ})_{\QQ}$ as follows.
Consider the central fonction $\alpha_{{\rm Con}}$ on $G$
defined by
$$\alpha_{{\rm Con}} (x) = 
\begin{cases}
1 & \text{if  $\langle x \rangle \in {\rm Con}$}\\
0 & \text{if  $\langle x \rangle \notin {\rm Con}$}.
\end{cases}
$$
Clearly 
$\alpha_{{\rm Con}}$ belongs to  $C (G, \QQ)$, hence,
by Proposition \ref{qval},
$\chi_{c} (C, \alpha_{{\rm Con}})$ belongs to
$K_{0}^{v} ({\rm Mot}_{k, \bar \QQ})_{\QQ}$. We
set
$\chi_{c} (C / A, {\rm Con}) := \chi_{c} (C, \alpha_{{\rm Con}})$.

Now let $X$ be a variety over $k$ and let  $\cA 
= \langle X, C_{i} / A_{i}, {\rm Con} (A_{i})
\, \vert \, i \in I \rangle$  be a Galois stratification of $X$.
We define the element
$$
\chi_{c} (\cA) := \sum_{i \in I} \chi_{c} (C_{i} / A_{i},
{\rm Con}(A_{i}))
$$
in 
$K_{0}^{v} ({\rm Mot}_{k, \bar \QQ})_{\QQ}$.

\begin{prop}\label{ind-is}Let $X$ be a variety over $k$ and let
$\cA$  be a Galois stratification of $X$.
The element $\chi_{c} (\cA)$ in $K_{0}^{v} ({\rm Mot}_{k, \bar \QQ})_{\QQ}$
depends only on 
the isomorphism class of the Galois stratification
$\cA$ of $X$.
\end{prop}

\begin{proof}Follows directly from the additivity of $\chi_c$ and Proposition
\ref{BN} (3).
\end{proof}

\begin{prop}\label{ext}Let
$k$ be a field of characteristic zero
which is the field of fractions of a normal domain $R$
of finite type over $\ZZ$.
Let $A$ be a normal integral variety
over $U := {\rm Spec} \, R$ and
consider  a Galois stratification $\cA 
= \langle A, C / A, {\rm Con} (A)\rangle$ 
consisting of a single colored Galois cover. Let $C' \rightarrow C$ be a
Galois cover such that the induced map
$C' \rightarrow A$ is a Galois cover. Consider  the family ${\rm Con}' (A)$ of subgroups
of $G (C'  / A)$ which 
are mapped onto
subgroups in ${\rm Con} (A)$ by the projection
$\varrho : G (C'  / A)  \rightarrow G (C  / A)$ and denote by
$\cA'$ the Galois stratification
$\cA' 
= \langle A, C' / A, {\rm Con}' (A)\rangle$.
\begin{enumerate}
\item[(1)]
For every   closed point $x$ of $U$ and every field $M$
containing $\FF_{x}$,
$
Z (\cA, x, M) = Z (\cA', x, M).
$
\item[(2)]We have
$$
\chi_{c} (\cA \otimes k) = \chi_{c} (\cA' \otimes k).
$$
\end{enumerate}
\end{prop}

\begin{proof}The first statement is clear and the second follows from
Proposition \ref{BN} (1), since $\alpha_{{\rm Con}' (A)} = 
\alpha_{{\rm Con} (A)} \circ \varrho$.
\end{proof}

\subsection{Frobenius action}\label{Frobenius action}We assume in this
subsection that
$k$ is a field of characteristic zero
which is the field of fractions of a normal domain $R$
of finite type over $\ZZ$.
We set again $U = {\rm Spec} \, R$.

Let us fix a prime number $\ell$ and denote by $G_{k}$ the absolute
Galois 
group
of $k$. We denote by
$K_{0} (\bar \QQ_{\ell}, G_{k})$ the Grothendieck 
group of the abelian category of finite dimensional $\bar 
\QQ_{\ell}$-vector spaces with continuous $G_{k}$-action. For every
closed point $x$ of $U$, we denote by ${\rm Frob}_{x}$ the 
geometric Frobenius automorphism over the field $\FF_{x}$. Taking the
trace of
${\rm Frob}_{x}$ on the invariants by inertia gives rise to a ring 
morphism 
$$
K_{0} (\bar \QQ_{\ell}, G_{k}) \otimes
\QQ \longrightarrow \bar \QQ_{\ell}.
$$
By composition with the morphism 
$$
{\rm \acute{E}t}_{\ell} : K_{0}^{v} ({\rm Mot}_{k, \bar \QQ}) \otimes \QQ
\longrightarrow K_{0} (\bar \QQ_{\ell}, G_{k}) \otimes
\QQ
$$
induced by {\'e}tale $\ell$-adic 
realization,
one defines a ring morphism
$$
{\rm Tr} \, {\rm Frob}_{x} : K_{0}^{v} ({\rm Mot}_{k, \bar \QQ}) \otimes \QQ
\longrightarrow \bar \QQ_{\ell}.
$$

\begin{prop}\label{frob}Let $X_{U}$ be a variety over
$U$ and let $\cA$ be a Galois stratification of
$X_{U}$. There exists a non zero element
$f$ in $R$ such that, for every closed point $x$ of $U_{f}$,
$${\rm Tr} \, {\rm Frob}_{x} (\chi_{c} (\cA \otimes k)) = {\rm card} \,
Z (\cA, x, \FF_{x}).$$
\end{prop}

\begin{proof}Follows directly from the next lemma.
\end{proof}

\begin{lem}
Let $A$ be a normal variety over $U$,
let $h : C \rightarrow A$
be a Galois cover with Galois group $G$ and let
${\rm Con}$ be a family of subgroups of
$G$ 
which is stable by conjugation under elements in $G$.
There exists a non zero element
$f$ in $R$ such that, for every closed point $x$ of $U_{f}$,
\begin{equation}\label{trace}{\rm Tr} \, {\rm Frob}_{x} (\chi_{c} (C / A
\otimes k,
{\rm Con}))
=
{\rm card} \, \Bigl\{a \in A (\FF_{x}) \Bigm \vert
{\rm Ar} (C / A, x, a) \subset
{\rm Con}
\Bigr\}.
\end{equation}
\end{lem}

\begin{proof}By its very
definition,  $\chi_{c} (C / A \otimes k, {\rm Con})$ is equal to
$\chi_{c} (C \otimes k, \alpha_{{\rm Con}})$.
As in the proof of 
Proposition \ref{qval} we may write $\alpha_{{\rm Con}}$
as a $\QQ$-linear combination
$$\alpha_{{\rm Con}} = \sum_{H} n_{H} {\rm Ind}^G_{H} 1_{H},$$
with $n_{H}$ in $\QQ$ and $H$ running over the set of cyclic
subgroups of $G$,
and it follows from Proposition \ref{BN} that
the left hand side of (\ref{trace}) is equal
to 
\begin{equation}\label{lhs}
\sum_{H} n_{H} \sum_{i} (-1)^{i} {\rm Tr} \, {\rm Frob}_{x}
(H^{i}_{c} ( C / H, \bar\QQ_{\ell})).
\end{equation}
The right hand side of (\ref{trace}) is equal
to $$
\sum_{a \in A (\FF_{x})} \alpha_{\rm Con} ({\rm Frob}_{a})
$$
with ${\rm Frob}_{a}$ the Frobenius automorphism at $a$
(up to conjugation)
and hence may be rewritten as
\begin{equation}\label{rhs}
\sum_{H} n_{H} {\rm card} \, \{a \in C / H (\FF_{x})\}.
\end{equation}
Now it follows
from Grothendieck's trace formula together with standard
construc\-ti\-bi\-li\-ty and base change theorems for $\ell$-adic cohomology
that there 
exists a non zero element
$f$ in $R$ such that, for every closed point $x$ of $U_{f}$,
(\ref{lhs}) is equal to (\ref{rhs}).
\end{proof}

In fact, the sets $Z (\cA, x, \FF_{x})$ completely determine the 
virtual
motive $\chi_{c} (\cA \otimes k)$ of a Galois 
stratification over $U$. More precisely:

\begin{theorem}\label{frobisom}Let $U = {\rm Spec} \, R$
be an affine scheme of finite 
type over $\ZZ$, integral, normal and of characteristic zero.
Let $X_{U}$ be a variety over
$U$ and 
let $\cA$ and $\cA'$
be Galois stratifications of $X_{U}$. Let $f$ be a non zero 
element of $R$. Assume that, for every closed point
$x$ of $U_{f}$, the equality
$Z (\cA, x, \FF_{x}) = Z (\cA', x, \FF_{x})$ holds. Then 
$$
\chi_{c} (\cA \otimes k)
=
\chi_{c} (\cA' \otimes k).
$$
\end{theorem}

\begin{proof}After refining $\cA$ and $\cA'$ one may assume
$I = I'$ and $A_{i} = A_{i}'$ for $i \in I$.
Hence it is sufficient to prove the following: let $A$
be an integral and  normal $U$-scheme of finite type and let
$h : C \rightarrow A$ and
$h' : C' \rightarrow A$ be Galois covers with Galois groups
respectively $G$ and $G'$ provided with
a family of subgroups ${\rm Con}$ (resp. ${\rm Con}'$) of
$G$ (resp. $G'$)
which is stable by conjugation under elements in $G$
(resp. $G'$); assume there exists a non zero element $f$ in $R$
such that for every closed point $x$ in $U_{f}$ and every
closed point $y$ of $A_{x}$,
${\rm Ar} (C / A, x, y) \subset {\rm Con}$
if and only 
${\rm Ar} (C' / A, x, y) \subset {\rm Con'}$,
then $$\chi_{c} (C \otimes k, \alpha_{{\rm Con}}) =
\chi_{c} (C' \otimes k, \alpha_{{\rm Con}'}).$$
But this follows from the more general Lemma \ref{chebo1}.
\end{proof}

\begin{lem}\label{chebo1}Let $A$
be an integral and  normal $U$-scheme of finite type and let
$h : C \rightarrow A$ and
$h' : C' \rightarrow A$ be Galois covers with Galois groups
respectively $G$ and $G'$. Take $\alpha \in \cC (G, \QQ)$
and
$\alpha' \in \cC (G', \QQ)$. Assume there exists a non zero $f$ in $R$
such that for every closed point $x$ in $U_{f}$ and every
closed point $y$ of $A_{x}$,
$$\alpha ({\rm Frob}_{y}) = 
\alpha' ({\rm Frob}_{y}).$$
Then $$\chi_{c} (C \otimes k, \alpha) =
\chi_{c} (C' \otimes k, \alpha').$$
\end{lem}

\begin{proof}Replacing $\alpha$ and $\alpha'$ by multiples, we may
assume that
$\alpha$ and $\alpha'$ are $\ZZ$-linear 
combinations of irreducible characters of $G$ and $G'$,
respectively.
By Proposition \ref{ext}, maybe after replacing $U$ by $U_{f}$
for some non zero $f$ in $R$,
we can assume $C = C'$ and $G = G'$ by going to a suitable commun
Galois
cover of $C$ and $C'$. By Chebotarev's Theorem, we have then
$\alpha =
\alpha'$,
and the result follows.
\end{proof}

More generally, the following result holds.

\begin{theorem}\label{superfrobisom}Let $U = {\rm Spec} \, R$
be an affine scheme of finite 
type over $\ZZ$, integral, normal and of characteristic zero.
Let $A$ and $B$ be varieties over
$U$ and 
let $\cA$ and $\cB$
be Galois stratifications of $A$ and $B$, respectively. 
Assume there exists a Galois stratification $\cZ$ of
$A \times B$ and 
a non zero 
element $f$ of $R$ such that, for every closed point
$x$ of $U_{f}$, $Z (\cZ, x, \FF_{x})$ is the graph of a bijection
between
$Z (\cA, x, \FF_{x})$ and $Z (\cB, x, \FF_{x})$.
Then 
$$
\chi_{c} (\cA \otimes k)
=
\chi_{c} (\cB \otimes k).
$$
\end{theorem}

\begin{proof}It is enough to prove that
$\chi_{c} (\cA \otimes k) = \chi_{c} (\cZ \otimes k)$.
Hence,
by additivity of $\chi_{c}$, one may assume that the
Galois stratification
$\cA$ consists of a single colored Galois cover $(C \rightarrow A, {\rm
Con} (A))$, with $A$ a normal variety and
${\rm Con} (A)$ non empty.
Denote by $Z$ the support of $\cZ$ and by
$p$
the restriction of the projection $A \times B \rightarrow A$
to $Z$. Write the
Galois stratification
$\cZ$ as
$$
\cZ
= \langle A \times B, W_{i} / Z_{i}, {\rm Con} (Z_{i})
\, \vert \, i \in I \rangle.
$$
It follows from Chebotarev's Theorem and the Lang-Weil estimate
\cite{Lang-Weil},
that the morphism $p : Z
\rightarrow A$ is generically finite and that
$\overline{p (Z)} = A$. 
Hence, maybe after performing
a finite partition  of $A$ into locally closed normal
subschemes, we may assume that,
for every $i$ in $I$, $Z_{i}$ is mapped by $p$ onto $A$
and that the restriction of $p$ to $Z_{i}$ is a finite {\'e}tale
morphism.
Furthermore, maybe after replacing $U$ by $U_{f}$ for a suitable non
zero $f$, we may,
by Proposition \ref{ext}, replace the $W_{i}$'s by suitable Galois
covers. It follows that, maybe after performing another
finite partition  of $A$ into locally closed normal
subschemes, we may assume
that the morphisms $W_{i} \rightarrow A$ obtained by
composition with $p$ are all Galois covers and 
that they are all equal to a same Galois cover $W \rightarrow A$.
Moreover we may assume that the cover $W \rightarrow A$ coincides with the
cover
$C \rightarrow A$.
By assumption there
exists
a non zero 
element $f$ of $R$ such that, for every closed point
$x$ of $U_{f}$, the map
$p$ 
induces a bijection between
$Z (\cZ, x, \FF_{x})$ and $Z (\cA, x, \FF_{x})$.
Hence, for such an $f$ and such an $x$,
a point $a$ belongs to $Z (\cA, x, \FF_{x})$
if and only if there exists $i$ in $I$ and $\sigma$ in
$G (W /A)$ such that
$\sigma \langle {\rm Frob}_a \rangle \sigma^{-1}$ belongs to
${\rm Con} (Z_i)$. Furthermore, since a given point $a$ in
$Z (\cA, x, \FF_{x})$ should be the image under $p$ of a unique point in
$Z (\cZ, x, \FF_{x})$,
one deduces from
Chebotarev's Theorem, that,
for $a$ in a subset of closed points of $A$ of Dirichlet density 1,
the element 
$i$ in $I$ and the class of $\sigma$ in $G (W / A) / G (W / Z_i)$
should both be unique. Since, for  such an $a$,
$$
\alpha_{{\rm Con}(A)} ({\rm Frob}_a)
=
\sum_{i \in I}
{\rm Ind}^{G (W / A)}_{G (W / Z_i)} \alpha_{{\rm Con} (Z_i)} ({\rm Frob}_a),
$$
it follows, again  by Chebotarev's Theorem, that
$$\alpha_{{\rm Con} (A)} =
\sum_{i \in I}
{\rm Ind}^{G (W / A)}_{G (W / Z_i)} \alpha_{{\rm Con} (Z_i)}.$$
By Proposition \ref{BN} (2), 
one deduces now that
$\chi_{c} (\cA \otimes k)
=
\chi_{c} (\cZ \otimes k)$.
\end{proof}

\subsection{Assigning virtual motives to formulas}\label{canmot}
We assume first that
$k$ is a field of characteristic zero which is of finite type over 
$\QQ$. We now associate a \rm{canonical}
element $\chi_{c} (\varphi)$ in
$K^{v}_{0} ({\rm Mot}_{k, \bar \QQ})_{\QQ}$ to every 
formula $\varphi$
in the first order language
of rings with coefficients in $k$ and free variables $X_{1}, \cdots,
X_{n}$.
The construction is as follows.
One may view $\varphi$ as a formula
in the first order language
of rings with coefficients in $R$, with 
$R$ a normal subring of $k$ which is of finite type over
$\ZZ$, and whose fraction field is $k$.
We set again $U = {\rm Spec} \, R$.
By Corollary \ref{2.5}  one may associate to $\varphi$ a Galois 
stratification 
$\cA$ of
$\AA^{n}_{U_{f}}$, for a suitable non zero $f$, such that,  for
every closed point $x$ of
$U_{f}$,
$$Z (\varphi, x, \FF_{x})
=
Z (\cA, x, \FF_{x}).$$ 
Clearly, the Galois stratification
$\cA \otimes k$ is not canonically associated to $\varphi$.
Nevertheless, it follows directly from
Theorem \ref{frobisom} that
the virtual motive 
$\chi_{c} (\cA \otimes k)$ in
$K^{v}_{0} ({\rm Mot}_{k, \bar \QQ})_{\QQ}$ is canonically attached to
$\varphi$. We shall denote it by $\chi_{c} (\varphi)$.

\medskip

Now we investigate the behaviour of 
$\chi_{c}$ under field extensions.

Let $k$ be a field of characteristic zero and
let $A$ be a normal integral variety over $k$.
Let $h : C \rightarrow A$
be a Galois cover with Galois group $G$ .
If  $k'$ is
a field containing $k$, the restriction of 
the morphism $h \otimes k': 
C \otimes k'\rightarrow A \otimes k'$ obtained by extension of scalars
to
a connected component $C'$ of 
of $C \otimes k'$ defines a Galois cover of $A \otimes k'$ with Galois
group $G'$, the decomposition subgroup of $G$ at the generic point of
$C'$, whose isomorphism class is independent of the choice of $C'$.
If the Galois cover is colored by a family ${\rm Con}$ of subgroups of
$G$ which is stable by conjugation under elements in $G$,
one considers the family ${\rm Con}'$ obtained by intersection with
$G'$. 
In this way, for $\cA$ a Galois
stratification
of a $k$-variety $X$, one defines a Galois
stratification $\cA \otimes k'$ of $X \otimes k'$.

\begin{lem}\label{preext}Let $k$ be a field of characteristic zero and
$\cA$ be a Galois stratification of a $k$-variety $X$.
Then $$\chi_c (\cA \otimes k') = \chi_c (\cA) \otimes k'.$$
\end{lem}

\begin{proof}It is enough to prove the following:
Let $h : C \rightarrow A$
be a Galois cover of $k$-varieties with Galois group $G$, $C'$
a connected component  of 
of $C \otimes k'$, $G'$ the decomposition subgroup of $G$ at the generic point of
$C'$, $\alpha$ a character of $G$, then
$\chi_c (C, \alpha) \otimes k' = \chi_c (C', {\rm Res}^G_{G'} \alpha)$.
Since $\chi_c (C, \alpha) \otimes k' = \chi_c (C \otimes k', \alpha)$,
this is equivalent to the equality
$\chi_c (C \otimes k', \alpha) = \chi_c (C', {\rm Res}^G_{G'} \alpha)$ which
follows from Proposition \ref{BN} (3).
\end{proof}

\begin{lem}\label{extension}Let $k \subset k'$ be fields
of finite type over $\QQ$. Let $\varphi$ be a formula
in the first order language
of rings with coefficients in $k$ and free variables $X_{1}, \cdots,
X_{n}$. We denote by $\varphi \otimes k'$ the same formula
considered as a formula with coefficients in $k'$.
Then
$$
\chi_{c} (\varphi \otimes k') = \chi_{c} (\varphi) \otimes k'.
$$
\end{lem}

\begin{proof}With the preceding notations take a Galois stratification
$\cA$ of
$\AA^{n}_{U_{f}}$, such that,  for
every closed point $x$ of
$U_{f}$,
$Z (\varphi, x, \FF_{x})
=
Z (\cA, x, \FF_{x})$. By Proposition \ref{2star},
$Z (\varphi, {\rm Spec} \, k, F)
=
Z (\cA,  {\rm Spec} \, k, F)$ for every pseudo-finite field $F$ containing
$k$. Hence, 
$Z (\varphi \otimes k', {\rm Spec}\,  k', F)
=
Z (\cA \otimes k',  {\rm Spec} \, k', F)$
for every pseudo-finite field $F$ containing
$k'$. Let 
$U' = {\rm Spec} \, R'$ with $R'$ the normalisation of $R$ in $k'$.
Again by
Proposition \ref{2star} we get that for some  non zero element $f'$ of $R'$,
$Z (\varphi \otimes k', x, \FF_{x})
=
Z (\cA \otimes U'_{f'}, x, \FF_{x})$, for
every closed point $x$ of
${U}'_{f'}$, with 
$\cA \otimes U'_{f'}$ the stratification obtained by base change.
Hence, by Theorem \ref{frobisom},
$\chi_{c} (\varphi \otimes k') = \chi_{c} (\cA \otimes k')$
and 
the result follows by
Lemma \ref{preext}.
\end{proof}

Let $k$ be any
field of characteristic zero. Let 
$\varphi$ be a formula in the first order language
of rings with coefficients in $k$. We may now
associate to $\varphi$
a \rm{canonical}
element $\chi_{c} (\varphi)$ in
$K^{v}_{0} ({\rm Mot}_{k, \bar \QQ})_{\QQ}$ as follows. Take a
subfield
$k_{0}$ of $k$ 
which is of finite type over 
$\QQ$ such that $\varphi$ may be viewed
as a formula $\varphi_{0}$ in the first order language
of rings with coefficients in $k_{0}$. By
Lemma \ref{extension}, $\chi_{c} (\varphi_{0}) \otimes k$ does not depend on
the choice of $k_{0}$, so we may set $\chi_{c} (\varphi)
:= \chi_{c} (\varphi_{0}) \otimes k$.

\begin{prop}\label{equiv}Let $k$ be a field of characteristic zero
and let $\varphi$ and $\varphi'$
be two formulas in the first order language
of rings with coefficients in $k$. If 
$\varphi \equiv \varphi'$ 
then
$$\chi_{c} (\varphi) = \chi_{c} (\varphi').$$
\end{prop}

\begin{proof}By Lemma \ref{rest} we may assume
that $k$ is of finite type over
$\QQ$, so we may apply Proposition
\ref{star} (2) to get a certain formula
$\psi$ and then, by Corollary \ref{2.5}, we may  replace the formula
$\psi$ by a Galois stratification and the result follows
from Theorem \ref{superfrobisom}.
\end{proof}

We shall denote by ${\rm Form}_{k}$ the set of formulas
in the first order language
of rings with coefficients in $k$.
If $\varphi$ is a formula in  ${\rm Form}_{k}$
with $n$ free variables
and
$\varphi'$ is a formula in  ${\rm Form}_{k}$ with
$n'$ free variables,
we shall denote by $\varphi \times \varphi'$ the formula 
with $n + n'$ free variables obtained by giving different names
to the free variables occuring in $\varphi$ and $\varphi'$
and taking the conjunction of $\varphi$
and $\varphi'$.

We now list a few properties of $\varphi \mapsto \chi_{c}(\varphi)$.

\begin{prop}\label{formulaire}The function $\chi_{c} :
{\rm Form}_{k} \rightarrow K^{v}_{0} ({\rm Mot}_{k, \bar \QQ})_{\QQ}$
satisfies the following properties.
\begin{enumerate}
\item[(1)]Let  $\varphi$ and $\varphi'$ be
formulas
with free variables $(X_{1}, \dots, X_{n})$. Then
$$\chi_{c} (\varphi \vee
\varphi') = \chi_{c} (\varphi) + \chi_{c} (\varphi') - \chi_{c} (\varphi \wedge
\varphi').$$
\item[(2)]Let $\varphi$ and $\varphi'$ be
formulas
with free variables $(X_{1}, \dots, X_{n})$ and $(X_{1}, \dots, X_{n'})
$,
respectively. Then
$$\chi_{c} (\varphi \times
\varphi') = \chi_{c} (\varphi)\chi_{c} (\varphi').$$
\item[(3)]If
$\varphi$
has $n$ free variables, then
$$\chi_{c} (\neg \varphi) = \LL^{n} - \chi_{c} (\varphi).$$
\item[(4)]$\chi_{c} (0 = 1) = 0$.
\end{enumerate}
\end{prop}

\begin{proof}To prove (1) first observe that
if $\varphi$ and $\varphi'$ are
formulas
with free variables $(X_{1}, \dots, X_{n})$, and $\chi_{c} (\varphi) =
\chi_{c} (\cA)$ and 
$\chi_{c} (\varphi') =
\chi_{c} (\cA')$ with $\cA$ and $\cA'$ Galois stratifications of
$\AA^{n}_{k}$
with
disjoint support then $\chi_{c} (\varphi \vee
\varphi') = \chi_{c} (\varphi) + \chi_{c} (\varphi')$.
Hence we have
$
\chi_{c} (\varphi \vee
\varphi') = \chi_{c} (\varphi \wedge
\varphi')  + \chi_{c} (\neg \varphi \wedge
\varphi')  + \chi_{c} (\varphi \wedge
\neg \varphi')
$,
$\chi_{c} (\varphi) = \chi_{c} (\varphi \wedge
\varphi')  +\chi_{c} (\varphi \wedge
\neg \varphi')$
and $\chi_{c} (\varphi') = \chi_{c} (\varphi \wedge
\varphi')  +\chi_{c} (\neg \varphi \wedge
\varphi')$, whence the statement follows.
(2) and (3) are proven by taking respectively
the product and the complement (in the obvious sense)
of the corresponding Galois
stratifications, while (4) is just evident.
\end{proof}

\subsection{New invariants of formulas}
Every ring morphism
$$K^{v}_{0} ({\rm
Mot}_{k, \bar \QQ})_{\QQ}
\longrightarrow R,$$ composed
with 
$\chi_{c} :
{\rm Form}_{k}  \rightarrow K_{0}^{v} ({\rm Mot}_{k, \bar \QQ})_{\QQ}$
will give rise to new invariants of formulas with coefficients in $k$.
Let us give
some examples. They all come from realization functors. Let $H^{\cdot}$ be a
cohomology theory on the category of smooth projective varieties over
$k$ with values in a field containing $\bar \QQ$. Then the realization
of a motive $(S, p, n)$ in 
${\rm Mot}_{k, \bar \QQ}$ is just $p(H^{\cdot}) \otimes H^{2}
(\PP^{1}_{k})^{\otimes n}$, with $p(H^{\cdot})$ the image of the projector $p$ acting 
on cohomology.
If one takes for $H^{\cdot}$ Betti or de Rham cohomology, taking the
alternating sum of the ranks of the cohomology groups
gives rise to the Euler characteristic ${\rm Eu} : K_{0}
({\rm Mot}_{k, \bar \QQ}) \rightarrow \ZZ$ and, after tensoring with
$\QQ$, to a morphism
${\rm Eu} : K_{0}
({\rm Mot}_{k, \bar \QQ}) \otimes \QQ \rightarrow \QQ$,
which restricts to a morphism
${\rm Eu} : K_{0}^{v}
({\rm Mot}_{k, \bar \QQ})_{\QQ} \rightarrow \QQ $.
Hence every formula
in ${\rm Form}_{k}$ has an Euler characteristic in $\QQ$.

By taking the Hodge realization, one associates to a motive
$M = (S, p, n)$ in 
${\rm Mot}_{k, \bar \QQ}$ its Hodge polynomial
$${\rm Hodge} (M) = \sum_{i, j} h^{i, j} (p(H^{\cdot}(S))) u^{i} v^{j} (uv)^{n}$$ in
$\ZZ [u, v, (uv)^{-1}]$, with $h^{i, j}$ the rank of the $(i, j)$-part
of $p(H^{\cdot})$. So we have a ring morphism
${\rm Hodge} : K_{0}
({\rm Mot}_{k, \bar \QQ}) \otimes \QQ \rightarrow \QQ[u, v, (uv)^{-1}]$,
which restricts to a morphism
${\rm Hodge} : K_{0}^{v}
({\rm Mot}_{k, \bar \QQ})_{\QQ} \rightarrow \QQ[u, v]$.
In this way one associates to  every 
formula in ${\rm Form}_{k}$ canonical Hodge numbers in
$\QQ$.

\medskip

We already considered in
\ref{Frobenius action} the Grothendieck 
group
$K_{0} (\bar \QQ_{\ell}, G_{k})$  of the abelian category of finite dimensional $\bar 
\QQ_{\ell}$-vector spaces with continuous $G_{k}$-action, and
the morphism 
$$
{\rm \acute{E}t}_{\ell} : K_{0}^{v} ({\rm Mot}_{k, \bar \QQ}) \otimes \QQ
\longrightarrow K_{0} (\bar \QQ_{\ell}, G_{k}) \otimes
\QQ
$$
induced by {\'e}tale $\ell$-adic 
realization. Hence to every  formula in
${\rm Form}_{k}$ is associated
a canonical virtual Galois representation with rational
coefficients in $K_{0} (\bar \QQ_{\ell}, G_{k}) \otimes
\QQ$.

\subsection{}
When $k$ is of finite type over $\QQ$, one has the following proposition.

\begin{prop}\label{frobfor}
Let $R$  be a normal domain  of finite type over
$\ZZ$, with fraction field $k$.
Let $\varphi$ be a formula
in the first order language
of rings with coefficients in $R$.
There exists a non zero element
$f$ in $R$ such that, for every closed point $x$ of ${\rm Spec} \, R_{f}$,
$${\rm Tr} \, {\rm Frob}_{x} (\chi_{c} (\varphi)) = {\rm card} \,
Z (\varphi, x, \FF_{x}).$$
\end{prop}

\begin{proof}Follows directly from
Proposition \ref{frob}.
\end{proof}

\begin{remark}The fact that the number of points of definable sets in
finite fields may be expressed as a $\QQ$-linear combination of
number of points of varieties goes back to \cite{Kiefe}.
For further results on the number of points of definable sets in
finite fields, see \cite{CDM} and \cite{F-H-J2}.
\end{remark}

\subsection{Grothendieck groups of first order theories}This subsection
is not used in the rest of the paper and its reading requires
some mild familiarity with the basic language of Model
Theory.

Let $\cL$ be a first order language and let $T$ be a theory in the
language 
$\cal L$.

We denote by $K_{0} (T)$ the quotient of the
free abelian group generated by
symbols $[\varphi]$ for $\varphi$ a formula in $\cL$
by the subgroup generated by the following
relations~
\begin{enumerate}
\item[(1)] If $\varphi$ is a formula in $\cL$ with free variables
$x = (x_{1}, \ldots, x_{n})$ and
$\varphi'$ is a formula in $\cL$ with free variables
$x' = (x'_{1}, \ldots, x'_{n'})$, then
$[\varphi] = [\varphi']$ if there exists a formula
$\psi$ in $\cL$, with free variables
$(x, x')$,
such that
$$
T \models [\forall x (\varphi (x) \rightarrow \exists! x' : (\varphi' (x')
\wedge \psi (x, x')))]
\wedge
[\forall x' (\varphi' (x') \rightarrow \exists! x : (\varphi (x)
\wedge \psi (x, x')))].$$
\item[(2)] \, $[\varphi \vee \varphi'] = [\varphi] + [\varphi'] - [\varphi
\wedge \varphi']$, for $\varphi$ and $\varphi'$  formulas in $\cL$.
\end{enumerate}
Furthermore one puts a ring structure on 
$K_{0} (T)$ by setting
\begin{enumerate}
\item[(3)] \, $[\varphi (x) ] \cdot [\varphi' (x')] =
[\varphi (x) \wedge  \varphi' (x')] $, if
$\varphi$ and $\varphi'$  are formulas in $\cL$ with disjoint free
variables $x$ and $x'$.
\end{enumerate}

For every interpretation of a theory $T_{1}$ in a theory
$T_{2}$ there is  a canonical morphism of rings
$K_{0} (T_{1}) \rightarrow K_{0} (T_{2})$, and this gives rise to a functor
from the category of theories in $\cL$, morphisms being given by
interpretation, to the category of commutative rings.

The previous constructions may now be rephrazed in the following way.
\begin{theorem}Let $k$ be a field of characteristic zero. Let $\cal L$
be the first order language of rings with coefficients in $k$ and let $T$ be
the theory of pseudo-finite fields containing $k$.
There exists a canonical morphism of rings 
$$
\chi_{c} : K_{0} (T)
\longrightarrow K_{0}^{v} ({\rm Mot}_{k, \bar \QQ})_{\QQ}$$ 
factorizing the morphism
$$
\chi_{c} : K_{0} ({\rm Sch}_{k})
\longrightarrow K_{0}^{v} ({\rm Mot}_{k, \bar \QQ})_{\QQ}.$$
\end{theorem}

\begin{proof}Indeed, this follows from
Proposition
\ref{equiv} and Proposition \ref{formulaire}.
\end{proof}

\begin{remark}If $k$ is a field and $T_{\rm ac}$ is 
the theory of algebraically closed fields containing $k$,
then $K_{0} (T_{\rm ac})$ is isomorphic to $K_{0} ({\rm Sch}_k)$. If
$T_{\RR}$ is the theory
of real closed fields in the language of ordered rings, 
then $K_{0} (T_{\RR})$ is isomorphic to $\ZZ$.
\end{remark}

\section{Definable subassignements for rings}\label{dsr}
\subsection{}\label{definable}
Let $h : \cC \rightarrow {\rm Sets}$ be a functor from a
category $\cC$ to the category of sets. We shall call
the data for each object $C$ of $\cC$ of a subset $h' (C)$ of
$h (C)$  a {\em subassignement} of $h$. The point in this definition
is that
$h'$ is not assumed to
be a subfunctor of $h$.

For $h'$ and $h''$ two subassignements
of $h$, we shall denote by $h' \cap h''$ and $h' \cup h''$, 
the subassignements
$ C \mapsto h' (C) \cap h'' (C)$ and $ C \mapsto h' (C) \cup h'' (C)$,
respectively. Similarly, we denote by $h \setminus h'$ the
subassignement $C \mapsto h (C) \setminus h' (C)$. We also write
$h \subset h'$ if $h (C) \subset h' (C)$ for every object $C$ of $\cC$.

Let $R$ 
be a ring.
We denote by ${\rm Field}_{R}$ the category
of fields which are $R$-algebras. For $X$ a variety over 
${\rm Spec} \, R$, we
consider the functor
$h_{X} : k \mapsto X (k)$ from ${\rm Field}_{R}$
to the category of sets.

Let $\varphi (X_{1}, \ldots, X_{m})$ be a formula 
in the first order
language of rings with coefficients in $R$ and free variables
$X_{1}, \ldots, X_{m}$. Let $\AA^{m}_{R}$ be the affine space ${\rm
Spec} \, R [X_{1}, \ldots, X_{m}]$.  For every field $k$ in
${\rm Field}_{R}$,
we denote by
$Z (\varphi, k) $ the subset of $k^{m} = \AA^{m}_{R} (k)$
defined by the formula $\varphi$.
This gives rise to a subassignement
$k \mapsto Z (\varphi, k) $
of the functor $h_{\AA^{m}_{R}}$.
We call such a subassignement a definable subassignement of 
$h_{\AA^{m}_{R}}$.
Let $X$ be a variety over ${\rm Spec} \, R$.
Assume first $X$ is affine and embedded as a closed subscheme in
$\AA^{m}_{R}$. We shall say a subassignement of $h_{X}$
is a definable subassignement if it is a definable 
subassignement of
$h_{\AA^{m}_{R}}$. Clearly, this definition is independent of
the choice of the embedding of $X$ in an affine space.
In general, a subassignement $h$ of $h_{X}$ will be
said to be definable if
there exists a finite cover
$(X_{i})_{i \in I}$ of $X$ by affine open subschemes and 
definable subassignements $h_{i}$ of $h_{X_{i}}$, for $i \in I$,
such that
$h  = \cup_{i \in I} h_{i}$.
When $X$ is affine this definition coincides with the previous one.

We shall denote by ${\rm Def}_{R} (X)$ the set of
definable subassignement of $h_{X}$. Clearly ${\rm Def}_{R} (X)$ is stable
by finite intersection and finite union and by taking complements.

\subsection{} Let now $k$ be a field of characteristic zero.
Assume first $X $ is a closed affine subvariety of $\AA^{m}_{k}$ and
let $h$ be a definable subassignement of $h_{X}$ associated to a formula
$\varphi$ with coefficients in $k$ and $m$ free variables. 
We set $\chi_{c} (h) = \chi_{c} (\varphi)$.
It follows
from Proposition \ref{equiv}
that $\chi_{c} (h)$ is independent of $\varphi$ and the embedding of
$X$. In general, when $X$ is a variety over $k$ and
$h$ is a definable subassignement of $h_{X}$, one takes a 
finite cover
$(X_{i})_{i \in I}$ of $X$ by affine open subschemes and 
definable subassignements $h_{i}$ of $h_{X_{i}}$ as in 
\ref{definable}
and one sets 
$$
\chi_{c} (h) = \sum_{\emptyset \not=J \subset I}
\chi_{c} (\cap_{i \in J} h_{X_{i}} \setminus \cup_{i \notin J} h_{X_{i}}),
$$
which is well defined by Proposition \ref{equiv} and
Proposition \ref{formulaire}. It also follows from
Proposition \ref{formulaire} that $\chi_{c}$ is additive,
{\it i.e.} that $\chi_{c} (h \cup h') = \chi_{c} (h) + \chi_{c} (h')
- \chi_{c} (h \cap h')$.

\begin{remark}
Of course we could also have defined $\chi_{c} (h)$ for
$h$ a definable subassignement of $h_{X}$ by directly
associating a Galois stratification to $h$, without considering
formulas.
\end{remark}

\subsection{}Let $R$  be  a normal domain of finite type over
$\ZZ$, with fraction field $k$. Let $X$ be a variety over $R$ and let
$h$ be a definable subassignement of $h_{X}$. We shall denote by
$h \otimes k$ the definable subassignement of $h_{X \otimes k}$ obtained by
extension of scalars.

The following result follows directly from Proposition \ref{frobfor}.

\begin{prop}\label{ppp}
Let $R$  be a normal domain  of finite type over
$\ZZ$, with fraction field $k$. Let $X$ be a variety over $R$ and let
$h$ be a definable subassignement of $h_{X}$.
There exists a non zero element
$f$ in $R$ such that, for every closed point $x$ of ${\rm Spec} \, R_{f}$,
$${\rm Tr} \, {\rm Frob}_{x} (\chi_{c} (h \otimes k)) = {\rm card} \,
h (\FF_{x}). \qed$$
\end{prop}

\subsection{}Let $R$ be a ring and 
$X$ be an  $R$-variety.
Let $h$ and $h'$ be definable subassignements
of $h_{X}$.
We shall write $h \approx h'$ if $h (k) = h' (k)$ for every
$k$ in ${\rm Field}_{R}$ which is a pseudo-finite field.
More generally, if $X$ and $X'$ are 
$R$-varieties, and $h$ and $h'$ are definable subassignements of
$h_{X}$ and $h_{X'}$, respectively, we write
$h \equiv h'$ if there exists a definable subassignement
$h''$ of $h_{X \otimes X'}$ such that, 
for every field
$k$ in ${\rm Field}_{R}$ which is a pseudo-finite field,
$h''(k)$ is the graph of a bijection between
$h (k)$ and $h' (k)$.
Clearly $\approx$ and $\equiv$ are equivalence relations.

Let $R$ be a normal domain  and let $X$ be a
variety
over $R$. Let $\cA$ be a Galois stratification of $X$. 
We can associate to $\cA$ a subassignement $h_{\cA}$ of
$h_{X}$ by defining, for $k$ in ${\rm Field}_{R}$,
as in \ref{Galois formulas},
$$
h_{\cA} (k) :=
\Bigl\{\aa = (a_{1}, \ldots, a_{n}) \in X (k)
\Bigm \vert {\rm Ar} (\aa) \subset
{\rm Con} (\cA)
\Bigr\}.
$$
We shall call such subassignements Galois subassignements
of $h_{X}$.

In this language Corollary \ref{pf2.5} can be reformulated as follows.

\begin{prop}Let $k$ be a field and let $X$ be a variety over $k$.
Let $h$ be a definable subassignement of $h_{X}$. Then there exists a Galois
subassignement $h'$ of $h_{X}$ such that $h \approx h'$.\qed
\end{prop}

\medskip

The following result follows directly from
Proposition \ref{equiv}.

\begin{prop}\label{hequiv}Let $k$ be a field and let $X$ and $X'$ be
$k$-varieties. If $h$ and $h'$ are definable subassignements of
$h_{X}$ and $h_{X'}$, respectively, such that $h \equiv h'$,
then $\chi_{c} (h) = \chi_{c} (h')$. \hfill \qed
\end{prop}

\subsection{}
Let $h$ be a definable subassignement
of $h_{X}$ for $X$ a variety over $k$. We say that 
$h$ is of dimension $\leq r$ if there exists a closed subvariety $S$ of
$X$ of dimension $\leq r$ such that $h $ is a subassignement of $h_{S}$.

\section{Definable subassignements  for power series rings}\label{dspr}

\subsection{Quantifier elimination for valued fields}
Let $K$ be a valued field, with valuation
${\rm ord} : K \rightarrow \Gamma \cup \{\infty \}$,
where $\Gamma$ is an ordered abelian group.
We denote by $\cO_{K}$ the valuation ring, by $P$ the valuation ideal,
by $U$ the group of units in $\cO_{K}$, by $\kappa$ the residue field, 
and by ${\rm Res} : \cO_{K} \rightarrow \kappa$ the canonical projection.
We assume that $K$ has an angular component map.
By this we mean a map ${\overline {\rm ac}} :
K \rightarrow \kappa$ such that
$\ac 0 = 0$, the restriction of $\ac$ to
$K^{\times}$ is multiplicative and 
the restriction of $\ac$ to $U$
coincides with the restriction of ${\rm Res}$.
From now on we fix that angular component map
${\overline {\rm ac}}$.

We consider 3-sorted first order languages\footnote{see,
{\it e.g.}, \cite{logic}
pp. 277---281 for more information concerning many-sorted first order
logic}
of the form 
$$
\cL = (\LL_{K}, \LL_{\kappa}, \LL_{\Gamma}, \ord, \ac),
$$
consisting of
\begin{enumerate}
\item[(i)] the language $\LL_{K} = \{+, -, \times, 0, 1\}$ of rings as 
valued
field sort,
\item[(ii)] the language $\LL_{\kappa} = \{+, -, \times, 0, 1\}$ of rings as 
residue field sort,
\item[(iii)] a language $\LL_{\Gamma}$, which is an extension of the 
language $\{+, 0, \infty, \leq\}$ of ordered abelian groups with an 
element $\infty$
on top, as the value sort,
\item[(iv)] a function symbol $\ord$ from the valued field
sort to
the value sort, which stands for the valuation,
\item[(v)]  a function symbol $\ac$ from the valued field sort to
the residue field sort, which stands for the angular component map.
\end{enumerate}

In the following we shall assume that $K$ is henselian and that $\kappa$ 
is of characteristic zero.
We consider $(K, \kappa, \Gamma \cup \{\infty\}, \ord, \ac)$ as
a structure
for the language $\cL$, the interpretations of symbols being the 
standard ones. By an  henselian
$\cL$-extension of $K$, we mean an
extension 
$(K', \kappa', \Gamma' \cup \{\infty\}, \ord', \ac')$ of the structure
$(K, \kappa, \Gamma \cup \{\infty\}, \ord, \ac)$ with respect to the language
$\cL$, with $K'$ a henselian valued field. By abuse of language we shall
say that $K'$ is a henselian $\cL$-extension of $K$.

We may now state the quantifier elimination Theorem of Pas 
\cite{P}.

\begin{theorem}\label{Pas}Let $K$ be a valued field
which satisfies the previous
conditions. For every  $\cal L$-formula $\varphi$ there exists 
an $\cal L$-formula $\varphi'$
without quantifiers over the valued field sort such that
$\varphi$
is equivalent
in $K'$ to $\varphi'$, for every
henselian $\cL$-extension $K'$ of $K$.
\end{theorem}

\begin{proof}This follows from Theorem 4.1 of \cite{P} together
with the observation at the begining of \S\kern .15em 3 of \cite{P}.
\end{proof}

In particular, when the value group is $\ZZ$, we shall use the language
$$
\cL_{\rm Pas} = (\LL_{K}, \LL_{\kappa}, \LL_{\rm PR \infty}, \ord, \ac),
$$
where $\LL_{\rm PR \infty} = \LL_{\rm PR} \cup \{\infty\}$ and $\LL_{\rm 
PR}$
is the Presburger language
$$
\LL_{\rm PR} = \{+, 0, 1, \leq\} \cup \{\equiv_{n} \, \vert \, n \in \NN, 
n > 1\},
$$
where $\equiv_{n}$ will be intepreted as ``congruent modulo $n$" in 
$\Gamma$.
We call a  subset of $\ZZ^{n}$ which is
definable in the language $\LL_{{\rm PR}}$
a Presburger subset of $\ZZ^{n}$. Similarly,
we call a function $\ZZ^{m} \rightarrow \ZZ^{r}$ a 
Presburger function if its 
graph
is definable in $\LL_{{\rm PR}}$.

\begin{cor}\label{prelq}Let $K$ be a valued field with value group
$\Gamma$
elementary equivalent\footnote{see, {\it e.g.}, \cite{F-J} for this notion} to $\ZZ$ in the language of ordered abelian
groups
and
satisfying the previous
conditions for
$\cL$ replaced by the language $\cal L_{\rm Pas}$.
For every  $\cal L_{\rm Pas}$-formula $\varphi$ there exists 
an $\cal L_{\rm Pas}$-formula $\varphi'$
without quantifiers over the valued field sort and the value sort
such that
$\varphi$
is equivalent
in $K'$ to $\varphi'$, for every
henselian $\cL_{\rm Pas}$-extension $K'$ of $K$
with value group
elementary equivalent to $\ZZ$.
\end{cor}

\begin{proof}Indeed, it follows from Theorem \ref{Pas},
since, by a classical result of Presburger \cite{Pr}, $\ZZ$ has
quantifier elimination in the language $\LL_{\rm PR}$.
\end{proof}

\subsection{}We assume from now
on that $K = k ((t))$, that
$\kappa = k$, with $k$
a field of characteristic zero,
and that $\ord$ and $\ac$ have their classical meaning for 
formal power series:
if $\varphi$ belongs to $k ((t))$, $\ord (\varphi)$ will denote the 
order in $t$ of $\varphi$ and $\ac (\varphi)$ the coefficient
of $t^{\ord (\varphi)}$ in $\varphi$, with the convention $\ac (0) = 0$.
In particular, the hypothesises of Corollary  \ref{prelq}
are satisfied in the language $\cL_{\rm Pas}$.

\medskip 
Let $R$ be a subring of $k$.
By {\emph {an $\cL_{\rm Pas}$-formula
with coefficients in $R$ in the valued field sort and in the residue
field sort}} we mean a formula in the language obtained from
$\cL_{\rm Pas}$ by adding, for every element of $R$,
a new symbol to denote it in the valued field sort and in the residue
field sort. Note that Corollary  \ref{prelq} remains valid for such formulas.
We shall consider $\cL_{\rm Pas}$-formulas
with coefficients in $R$ in the valued field sort and in the residue
field sort, free
variables $x_{1}, \ldots, x_{m}$ running
over the valued field sort and no free
variables running over
the residue field or the value sort. We shall call such formulas
\textit{formulas on $R [[t]]^m$}.
The reason for that denomination is that later
we shall view the free variables $x_{1}, \ldots, x_{m}$
as running over $R[[t]]$. 
More generally, an $\cL_{\rm Pas}$-formula
with coefficients in $R$ in the valued field sort and in the residue
field sort, free
variables $x_{1}, \ldots, x_{m}$ running
over the valued field sort, no free
variables running over
the residue field sort and $r$ free variables running over the value sort,
gives rise by specialization of the value sort variables to 
$\ZZ^r$ (resp. $\NN^{r}$ or $\NN \cup \{\infty\}$ when $r = 1$)
to what we shall call a \textit{formula on $R [[t]]^m$
depending on parameters
in} $\ZZ^r$ (resp. $\NN^{r}$ or $\NN \cup \{\infty\}$).

\medskip
We shall deduce the
following statement  of Ax/Ax-Kochen-Er{\v s}ov type 
from
the Theorem of Pas.

\begin{prop}\label{axax}Let $R$ be a
normal domain 
of finite type over $\ZZ$ with field of fractions $k$.
Let $\sigma$ be a sentence
in the language $\cL_{{\rm Pas}}$ with coefficients in $R$
in the valued field sort and in the residue
field sort. The following
statements are equivalent:
\begin{enumerate}
\item[(1)] The sentence
$\sigma$ is true in $F [[t]]$ for every pseudo-finite field $F$
containing
$k$.    
\item[(2)]
There exists $f$ in $R \setminus \{0\}$
such that, for every closed point $x$ in ${\rm Spec} \, R_{f}$,
the sentence $\sigma$ is true in $\FF_{x}[[t]]$.
\end{enumerate}
If, furthermore, $k$ is a finite extension of $\QQ$, the previous statements
are also equivalent to the following:
\begin{enumerate}
\item[(3)]There exists $f$ in $R \setminus \{0\}$, multiple of the
discriminant of $k / \QQ$,
such that, for every closed point $x$ in ${\rm Spec} \, R_{f}$,
the sentence $\sigma$ is true in $k_{x}$,
\end{enumerate}
where $k_{x}$ denotes the completion of $k$ at $x$. Remark that, 
the extension $k / \QQ$ being non ramified at $x$,
the field $k_{x}$ admits a canonical uniformizing parameter, hence
also
a canonical
angular component map.
\end{prop}

\begin{proof} Let us first prove the equivalence of (1) and (2).
By Corollary \ref{prelq}, there exists
an $\cal L_{\rm Pas}$-sentence $\sigma'$
without quantifiers over the valued field sort and the value sort
such that
$\sigma$
is equivalent
in $K'$ to $\sigma'$, for every
henselian $\cL_{\rm Pas}$-extension $K'$ of $k [[t]]$. Hence
there exists $f$ in $R \setminus \{0\}$
such that, for every closed point $x$ in ${\rm Spec} \, R_{f}$,
the sentence $\sigma$ is equivalent to  $\sigma'$
in $\FF_{x}[[t]]$.
Indeed, if this would not be the case,
a suitable ultraproduct of
the fields $\FF_{x}[[t]]$ would yield a henselian
$\cL_{\rm Pas}$-extension $K'$ of $k [[t]]$ in which $\sigma$
would not be  equivalent
to $\sigma'$. Hence we may assume $\sigma$ is a sentence
({\it i.e} is without free variables) and has 
quantifiers only over the residue field sort, in which case the
result follows from
Proposition \ref{2star}. The proof of the equivalence of (1) and (3)
is completely similar.
\end{proof}

\subsection{}We shall call a formula on $R [[t]]^m$ a \textit{special formula on $R [[t]]^m$}
if it is obtained by 
repeated application of
conjunction and negation from formulas of the form
\begin{gather}
 \ord f_1 (x_1, \ldots, x_m) \geq
\ord f_2 (x_1, \ldots, x_m) + a\label{aaa}\\
\ord f_1 (x_1, \ldots, x_m) \equiv
a \mod b\label{bbb}
\intertext{and}
\vartheta (\ac (f_1 (x_1, \ldots, x_m)),
\dots,
\ac (f_{m'} (x_1, \ldots, x_m))), \label{ccc}
\end{gather}
where the $f_i$ 
are polynomials with coefficients in  $R$, 
$a$ and $b$ are in $\ZZ$, and  $\vartheta$ is a formula in $m'$ free variables
in the first order language of rings with coefficients 
in $R$.

Replacing $a$ in the above formulas by $L (a_{1}, \dots, a_{r})$
with $L$ a polynomial with coefficients in $\ZZ$ and degree $\leq 1$, 
one gets the definition of   
\textit{special formulas on $R [[t]]^m$
depending on parameters}
in $\ZZ^r$ (resp. $\NN^{r}$  or $\NN \cup \{\infty\}$).

\subsection{}
Let $k$ be a field.
For $X$ a variety over $k$, we will denote by $\cL (X)$ the scheme of 
germs of
arcs on $X$. It is a scheme over $k$ and, for every field extension
$k \subset K$, there is a natural bijection,
$$\cL (X) (K) \simeq {\rm Mor}_{k-{\rm schemes}} (\Spec K [[t]], X),
$$ 
between the set of $K$-rational points of $\cL (X)$
and the set of 
germs of arcs with coefficients in  $K$ on $X$.
We will call the
$K$-rational points of $\cL (X)$, for $K$
a field extension of $k$, arcs on $X$, and $\varphi (0)$ will be called 
the origin of the arc $\varphi$.
More precisely, the
scheme $\cL (X)$ is defined as the projective limit,
$$
\cL (X) := \varprojlim \cL_{n} (X),
$$
in the category of $k$-schemes of the schemes
$\cL_{n}(X)$ representing the functor,
$$R \mapsto {\rm Mor}_{k-{\rm schemes}} (\Spec R [t] / t^{n+1} R[t], X),$$
defined on the category of $k$-algebras. (The existence of $\cL_{n}(X)$
is well known, cf. \cite{Arcs}, and the projective limit exists since 
the transition morphisms are affine.)
We shall denote by $\pi_{n}$ the canonical morphism corresponding to 
truncation of arcs,
$$
\pi_{n} : \cL (X) \longrightarrow \cL_{n} (X).
$$
The schemes $\cL (X)$ and $\cL_{n} (X)$ will always be
considered to have their reduced structure.

\subsection{}Let $R$ be a ring
and let $X$ be a variety over ${\rm Spec} \, R$. 
We consider the functor
$h_{\cL (X)} :  k \mapsto X (k[[t]])$ from ${\rm Field}_{R}$
to the category of sets.

Let $\varphi$ be a formula on $R [[t]]^m$.
For every field $k$ in ${\rm Field}_{R}$, denote by 
$Z (\varphi, k[[t]]) $ the subset of $k[[t]]^{m} = \AA^{m}_{R} (k [[t]])$
defined by the formula $\varphi$.
This defines a subassignement
$k \mapsto Z (\varphi, k[[t]]) $
of the functor $h_{\cL (\AA^{m}_{R})}$.
We call such a subassignement a definable subassignement of 
$h_{\cL (\AA^{m}_{R})}$.
We now proceed in a
similar way as in \ref{definable}
to define  definable subassignements of $h_{\cL (X)}$, for
$X$ a variety over ${\rm Spec} \, R$.
When $X$ is affine and embedded as a closed subscheme in
$\AA^{m}_{R}$, we  shall say a subassignement of $h_{\cL (X)}$
is a definable subassignement if it is a definable subassignement of
$h_{\cL (\AA^{m}_{R})}$. 
In general a subassignement $h$ of $h_{\cL (X)}$ will be
said to be definable if
there exists a finite cover
$(X_{i})_{i \in I}$ of $X$ by affine open subschemes and 
definable subassignements $h_{i}$ of $h_{\cL (X_{i})}$
such that,
$h  = \cup_{i \in I} h_{i}$.
Similarly, one defines definable subassignements $h_{a}$ of $h_{\cL (X)}$
depending on parameters $a = (a_{1}, \dots, a_{r})$ in
$\ZZ^r$, in $\NN^{r}$, or in $\NN \cup \{\infty\}$.

We shall denote by ${\rm Def}_{R} (\cL (X))$ the set of
definable subassignements of $h_{\cL (X)}$. Clearly ${\rm Def}_{R} (\cL
(X))$ is stable
by finite intersection and finite union and by taking complements.

We shall use the symbol $\buildrel {\cdot}
\over {\cup}$ to denote finite or infinite union of pairwise
disjoint subassignements of $h_{\cL (X)}$.

For a ring morphism $R \rightarrow R'$ and
$h$ a functor or a subassignement
from ${\rm Field}_{R}$ to the category of sets,
we denote by $h \otimes R'$ the restriction of $h$ to
${\rm Field}_{R'}$.

\medskip
The following statement is a direct consequence of
Corollary
\ref{prelq}.

\begin{prop}\label{Pas6}Let $k$ be a field of characteristic zero.
Let $\varphi$ be a formula (resp. a formula depending on parameters) on $k[[t]]^{m}$.
Then there exists a special formula (resp. a special formula depending on
parameters) $\varphi'$ on $k[[t]]^{m}$
such that $\varphi$ and $\varphi'$ define the same definable 
subassignement (resp. the same definable 
subassignement depending on parameters) of $h_{\cL
(\AA^{m}_{k})}$. \hfill \qed
\end{prop}

\subsection{}Let $k$ be a field and 
$X$ be a $k$-variety.
Let $h$  be a definable subassignement of $h_{\cL (X)}$. By a definable
partition of $h$, with parameters in $\ZZ$, $\NN$, or
in $\NN \cup \{\infty\}$,
we mean  the data of definable subassignements
$h_{i}$, depending on the
parameter $i \in \ZZ$, $\NN$, or in $\NN \cup \{\infty\}$, which are 
pairwise disjoint and such that $h = {\bigcup_{i}} h_{i}$.
Similarly, if $h_{n}$ already depends on a parameter $n$,
a definable partition of the subassignements $h_{n}$,
with parameters in $\ZZ$, $\NN$, or in $\NN \cup \{\infty\}$, will be
the data of definable subassignements $h_{n, i}$, depending on the
parameter $(n, i)$,
such that, for each $n$, the subassignements  $h_{n, i}$
are 
pairwise disjoint and 
$h_{n} = 
{\bigcup_{i}} h_{n, i}$.

\subsection{Truncation of definable subassignements}Let $k$ be a field
of characteristic zero
and let $X $ be a variety over $k$.
For $n$ in $\NN$, we have a canonical truncation 
morphism $\pi_{n} : \cL (X) \rightarrow \cL_{n} (X)$.
Hence if $h$ is a subassignement of $h_{\cL(X)}$ 
(resp. of $h_{\cL_{n} (X)}$)
we may consider
$\pi_{n} (h) : K \mapsto \pi_{n} (h (K))$
(resp. $\pi_{n}^{-1} (h) : K \mapsto \pi_{n}^{-1} (h (K))$)
which is a subassignement of $h_{\cL_{n} (X)}$
(resp. of $h_{\cL (X)}$).

\begin{prop}\label{speci}Let $h$ be  a definable
subassignement of $h_{\cL(X)}$. Then, for every 
$n$ in $\NN$, 
$\pi_{n} (h)$ is a definable
subassignement of $h_{\cL_{n} (X)}$ in the sense of \S\kern .15em \ref{dsr}
and
$\pi_{n}^{-1} \pi_{n} (h)$ is a definable
subassignement of $h_{\cL (X)}$.
\end{prop}

\begin{proof}One may assume $X = \AA^{m}_{k}$
and $h$ is associated to a formula
$\varphi$ on $k[[t]]^m$. Then $\pi_{n}^{-1} \pi_{n} (h)$
is associated to the formula
$$
\exists \, y_{1} \dots \exists \, y_{m} \,
\Bigl ((\ord (x_{1} - y_{1}) > n) \wedge
\dots \wedge
(\ord (x_{m} - y_{m}) > n) \wedge \varphi (y_{1}, \dots, y_{m}) \Bigr)
$$
and it follows from Corollary \ref{prelq} that $\pi_{n}^{-1} \pi_{n}
(h)$
is a definable subassignement of $h_{\cL (X)}$ associated to 
a special formula 
on $k[[t]]^m$ - without quantifiers in the valued field sort - say $\psi$.
One can formally write
$x_{i} = \sum_{j \geq 0} a^{(i)}_{j} t^j$ and then, by expanding
the variables $x_{i}$ into series in
the formulas of type (\ref{aaa}),(\ref{bbb}) and (\ref{ccc})
appearing in $\psi$, one 
obtains a infinite set of conditions  in the variables $a^{(i)}_{j}$. 
If one substitutes in these conditions  $a^{(i)}_{j} = 0$
whenever $i > n$, then
only a finite number of conditions which
involve  only coefficients
$a^{(i)}_{j}$ with $j \leq n$ remain.
In this way one obtains from the formula $\psi$ a formula
$\tau_{n} (\psi)$ with coefficients in $k$ and   free  variables
$a^{(i)}_{j}$, $j \leq n$ and $\pi_{n} (h)$ is the  definable
subassignement of $h_{\cL_{n} (X)}$ associated to $\tau_{n} (\psi)$.
\end{proof}

\subsection{Stable definable subassignements}Let $k$ be a field
of characteristic zero
and let $X $ be a variety over $k$.
Let $h$ be a definable subassignement of $h_{\cL (X)}$.
We say $h$ is weakly stable at level $n$ if
$h = \pi^{-1}_{n} \pi_{n} (h)$ and that
$h$ is weakly stable if it is weakly stable at some level.

\begin{lem}\label{2.4}Let  $h$ 
and $h_{i}$, $i \in \NN$, be weakly stable
definable subassignements of $h_{\cL (X)}$.
Assume that
$$
h = \bigcup_{i \in \NN} h_{i}.
$$
Then there exist a natural number $n$ such that
$$
h = \bigcup_{i \leq n} h_{i}.
$$
\end{lem}

\begin{proof}The  proof of Lemma 2.4 in
\cite{Arcs} using ultraproducts
may be directly adapted to the present situation.
Indeed, we may assume $X$ is a closed subvariety
of $\AA^m_k$, and it is enough to prove that
if $k_{i}$, $i \in \NN$, are weakly stable
definable subassignements of $h_{\cL (X)}$ such that, for every finite
subset $\Sigma$ of $\NN$, $\bigcap_{i \in \Sigma} k_i$ is not the
empty subassignement,  {\it i.e.} for some field $K_{\Sigma}$,
$\bigcap_{i \in \Sigma} k_i (K_{\Sigma})$ is not empty,
then $\bigcap_{i \in \NN} k_i$ is not the
empty subassignement.
Since every weakly stable
definable subassignement of $h_{\cL (X)}$
may be defined  by an infinite conjunction of formulas
in the language of rings with coefficients in $k$ involving each
only a
finite
number of coefficients of the power series $x_{i}$, $1 \leq i \leq
m$, cf. the proof of Proposition \ref{speci}, it follows that
$\bigcap_{i \in \NN} k_i (K^{\ast})$ is not empty for
$K^{\ast}$
the ultraproduct of the fields $K_{\Sigma}$
with respect to a suitable ultrafilter.
\end{proof}

\bigskip
Let $\pi : X \rightarrow Y$ be a morphism of
algebraic varieties over $k$
and let
$h$ and $h'$ be definable subassignements of
$h_{X}$ and $h_{Y}$, respectively. Assume $\pi (h) \subset h'$.
We say the $\pi$ induces a  {\it piecewise trivial fibration}
$h \rightarrow h'$ 
{\it with fiber} a $k$-variety 
$F$, if there exists a finite family
of locally closed subsets $S_{i}$, $i \in I$, of $Y$, such that 
$\pi^{- 1} (S_{i})$ is locally closed in $X$,
with $h' \subset \cup_{i \in I} h_{S_{i}}$, such that
there is an isomorphism
$\pi^{-1} (S_{i}) \simeq S_{i} \times F$,
with $\pi$
corresponding under the isomorphism to the projection 
$S_{i} \times F \rightarrow S_{i}$, inducing, for every $L$ in
${\rm Field}_{k}$, a bijection between
$(h \cap h_{\pi^{-1} (S_{i})}) (L)$ and
$(h' \cap h_{S_{i}}) (L) \times h_{F} (L)$.

\bigskip
Let $X$ be an algebraic variety over $k$ of 
dimension
$d \geq 0$
and let $h$ be a definable subassignement of $h_{\cL (X)}$. We say $h$
is
{\it stable at level} $n \in \NN$, if $h$ is weakly
stable at level $n$ and the canonical morphism
$\cL_{m + 1} (X) \rightarrow \cL_{m} (X)$ induces
a piecewise trivial fibration
$\pi_{m + 1} (h) \rightarrow \pi_{m} (h)$
with fiber
$\AA^{d}_{k}$ for every $m \geq n$.

We say $h$ is {\it stable} if it stable at some level $n$. 
The set ${\rm Def}_{k} (\cL
(X))_{st}$
of stable definable subassignements of
$h_{\cL (X)}$ is stable by 
taking finite intersections and finite unions.

\begin{remark}\label{4.1}
If $h$ is a weakly stable subassignement of $h_{\cL (X)}$ and
$h \cap h_{\cL (X_{0})} = \emptyset$,
with $X_{0}$ the union of the singular locus of $X$ and its irreducible
components of dimension $<d$, then
$h$ is stable, as follows from Lemma \ref{2.4} and Lemma 4.1 of \cite{Arcs}. 
In particular, if $X$ is smooth of pure dimension $d$,
every weakly stable subassignement of
$h_{\cL (X)}$ is stable.
\end{remark}

\begin{lem}\label{triviallemma}If $h$ is stable at level $n$, then,
for every $n' \geq n$,
$$
\chi_{c} (\pi_{n'} (h)) \, \LL^{ - (n' + 1) d} = 
\chi_{c} (\pi_{n} (h)) \, \LL^{ - (n + 1) d}.
$$
\end{lem}

\begin{proof}Clear.
\end{proof}

\subsection{}Let $R$ be a ring and 
$X$ be an $R$-variety.
Let $h$ and $h'$ be subassignements of $h_{\cL (X)}$.
We shall write $h \approx h'$ if $h (k) = h' (k)$ for every
$k$ in ${\rm Field}_{R}$ which is a pseudo-finite field.
If $h \approx h'$, then $\pi_{n} (h) \approx \pi_{n} (h')$ for every $n$
in $\NN$.

\section{Arithmetic motivic integration on definable sets}\label{amids}
\subsection{}Let $k$ be a field of characteristic 0 and let $X$ be a
variety over $k$ of  dimension $d$.
By Lemma \ref{triviallemma} there exists a unique function
$$
\tilde \nu : {\rm Def}_{k} (\cL (X))_{st}
\longrightarrow
K_{0}^{v} ({\rm Mot}_{k, \bar \QQ})_{{\rm loc}, \QQ}
$$
such that 
$$
\tilde \nu (h) = \chi_{c} (\pi_{n} (h)) \, \LL^{ - (n + 1) m}
$$
if $h$ is stable of level $n$.

\begin{prop}\label{add}Let $h$ and $h'$ be
in ${\rm Def}_{k} (\cL (X))_{st}$.
Then
\begin{enumerate}
\item[(1)]
$\tilde \nu (h \cup h') =
\tilde \nu (h) + \tilde \nu (h') -
\tilde \nu (h \cap h')
$.
\item[(2)]If $h \approx h'$,
then $ \tilde \nu (h) =  \tilde \nu (h')$.
\end{enumerate}
\end{prop}

\begin{proof}The first assertion follows directly from
Proposition \ref{formulaire}.
The second assertion is a consequence of Proposition \ref{hequiv},
together with the fact that $\pi_{n} (h) \approx \pi_{n} (h')$
if $h \approx h'$.
\end{proof}

Now let $h$ be a
stable definable subassignement of 
$h_{\cL (X)}$ and consider a definable
partition of $h$ with parameters in 
$\ZZ$ such that each $h_{n}$ is stable. Then, by Lemma \ref{2.4},
$\tilde \nu (h_{n}) = 0$ for $\vert n \vert \gg 0$ and the sum
$\sum_{n \in \ZZ} \LL^{- n} \tilde \nu (h_{n})$ is finite. 
We denote that sum by $\int \LL^{- n} h_{n} d \tilde \nu$.

\bigskip

In general, for $h$ a
definable subassignement of 
$h_{\cL (X)}$
which is maybe not stable, one has to use a limit
process to define a ``motivic measure'' of $h$.
To achieve this aim, one needs the following lemma.

\begin{lem}\label{cuting}Let $X$ be a
variety over $k$ of dimension $d$ and let
$h$ be a
definable subassignement of 
$h_{\cL (X)}$. Then there exists definable subassignements
$k_{i}$,  stable at level $n_{i}$, depending on the parameter $i \in \NN$,
and a closed subvariety $S$ of dimension $ < d$ of $X$
such that 
$$
h =  \quad \buildrel {\cdot} \over {\bigcup_{i}} 
k_{i} \buildrel {\cdot} \over {\cup} (h \cap h_{\cL (S)}),
\quad
\lim_{i \rightarrow \infty} ({\rm dim} \, \pi_{n_{i}}(k_{i}) 
- (n_{i} + 1) d) = - \infty,$$ and such
that the denominators of the elements
$\chi_{c} (\pi_{n_{i}}(k_{i}) ) $ in
$K_{0}^{v} ({\rm Mot}_{k, \bar \QQ})_{\QQ}$, $i \in \NN$, are bounded.
Furthermore, if $h_{n}$, $n \in \ZZ$,
is a definable
partition of $h$, it is possible to choose the
$k_{i}$'s in such a way that,
for every $n$, $h_{n}$ is contained in some
$k_{i}$.
\end{lem}

\begin{proof}One reduces first
to the case when
$X$ is affine irreducible with a closed immersion $X \rightarrow \AA^{m}_{k}$,
and then, by Proposition \ref{Pas6}, we may assume 
$h$ is associated to a special formula $\varphi$ on $k [[t]]^{m}$
which is 
obtained by repeated application of
conjunction and negation from formulas of the form
(\ref{aaa}),  (\ref{bbb}) and  (\ref{ccc}). Choose a nonzero regular
function
$g$ on $X$ which vanishes on the singular locus of $X$ and 
let $f$ denote the product of $g$ and all the polynomials
$f_{i}$, assumed to be non zero, occuring in (\ref{aaa}),  (\ref{bbb})
and  (\ref{ccc}).
Now set
$\psi_{i} = \varphi \wedge \cap ({\rm ord}_{t} f = i)$
and define $S$ as the locus of $f = 0$.
The definable subassignements $k_{i}$ defined by $\psi_{i}$
are stable by Remark \ref{4.1}
and the statement on
dimension
follows from Lemma 4.4 of \cite{Arcs}. We still have to check that
the denominators of
$\chi_{c} (\pi_{n_{i}}(k_{i}) ) $ in
$K_{0}^{v} ({\rm Mot}_{k, \bar \QQ})_{\QQ}$ are bounded.
For that it is enough to
know that one can bound uniformly the degree of the coverings
in the Galois
stratifications associated by quantifier elimination (Corollary
\ref{pf2.5}) to 
$\pi_{n_{i}}(k_{i})$. But consider  the Galois stratifications
associated by Corollary  \ref{pf2.5} to the formulas 
$\vartheta$ occuring in (\ref{ccc}) and let $d$ the maximum of the
degrees
of the coverings appearing in these Galois stratifications.
The Galois
stratifications associated to 
$\pi_{n_{i}}(\psi_{i})$ may be expressed
in terms of the former ones, and the integer $d$ is
still a bound for the
degrees
of the coverings. 
For definable partitions, the construction of the $\psi_{i}$'s
is done in exactly in the same way.
\end{proof}

\begin{def-theorem}\label{ami}There exists a unique
mapping
$$\nu : {\rm Def}_{k} (\cL (X))
\longrightarrow
\widehat K_{0}^{v} ({\rm Mot}_{k,
\bar \QQ})_{\QQ}$$
satisfying the following properties.
\begin{enumerate}
\item[(1)]If  $h$ is a stable definable subassignement of $h_{\cL (X)}$,
then
$\nu (h)$
is equal to the image
of
$\tilde \nu (h)$ in $\widehat K_{0}^{v} ({\rm Mot}_{k,
\bar \QQ})_{\QQ}$.
\item[(2)]If $h$ and $h'$ are definable subassignement of $h_{\cL (X)}$, then
$$\nu (h \cup h') =
 \nu (h) +  \nu (h') -
 \nu (h \cap h').
$$
\item[(3)]If $h \approx h'$,
then $ \nu (h) =   \nu (h')$.
\item[(4)]Let  $h$ be a definable subassignement of $h_{\cL (X)}$.
If ${\rm dim} \, h \leq d - 1$, then
$\nu (h) = 0$.
\item[(5)]Let $h_{n}$ be a definable partition of a definable subassignement
$h$
with parameter $n \in \NN$.
Then the series $\sum_{n \in \NN} \nu (h_{n})$
is convergent in
$\widehat K_{0}^{v} ({\rm Mot}_{k, \bar \QQ})_{\QQ}$
and 
$$
\nu (h) = \sum_{n \in \NN} \nu (h_{n}).
$$
\item[(6)]Let  $h$ and $h'$ be definable subassignements of $h_{\cL (X)}$.
Assume
$h \subset h'$.
If $ \nu (h') $ belongs to $F^{e} \widehat K_{0}^{v} ({\rm Mot}_{k, \bar
\QQ})_{\QQ}$, then $\nu (h)$
also belongs to $F^{e} \widehat K_{0}^{v} ({\rm Mot}_{k, \bar
\QQ})_{\QQ}$.
\end{enumerate}
We call $\nu (h)$ the arithmetic motivic volume
of $h$.
\end{def-theorem}

\begin{proof}The proof is just the same as the proof of
Definition-Proposition 3.2 in \cite{Arcs}, if one replace
Lemma 2.4 and Lemma 3.1 of loc. cit. by Lemma \ref{2.4} and
Lemma \ref{cuting}.
\end{proof}

Let $h_{n}$ be a definable partition of a definable subassignement $h$
with parameter $n \in \ZZ$.
We say that $\LL^{-n} h_{n}$ is integrable
if the series 
$$
\int \LL^{- n} h_{n} d\nu:= \sum_{n \in \ZZ} \LL^{- n} \nu (h_{n})
$$
converges in $\widehat K_{0}^{v} ({\rm Mot}_{k, \bar \QQ})_{\QQ}$.
It follows from \ref{ami} (6) that if $h_{n}$ is
a definable partition of a formula $h$
with parameter $n \in \NN$ (or in $\NN \cup \{\infty\}$ with the convention
$\LL^{- \infty} \nu (h_{\infty}) = 0$),
then $\LL^{-n} h_{n}$ is integrable.

\medskip

The following result is the analogue of
Theorem 7.1 of \cite{Arcs} in the present context.

\begin{theorem}Let $X$ be a variety over $k$ of dimension $d$.
Let  $h$ be a definable subassignement of $h_{\cL (X)}$. Then 
$$
\lim_{n \rightarrow \infty}
\chi_{c} (\pi_{n} (h)) \, \LL^{-(n +
1) d}
= 
\nu (h)
$$ in
$\widehat K_{0}^{v} ({\rm Mot}_{k, \bar \QQ})_{\QQ}$.
\end{theorem}

\begin{proof}Again the  proof is essentially  the same as the proof of
Theorem 7.1 in \cite{Arcs}, if one replace
Lemma 2.4 and Lemma 3.1 of loc. cit. by Lemma \ref{2.4} and
Lemma \ref{cuting}.
\end{proof}

\subsection{Change of variable formula}
Let $X$ be an algebraic variety over $k$ of dimension
$d$, and let $\cI$ be a coherent sheaf of ideals on $X$.
We denote by 
${\rm ord}_t \cI$
the function
${\rm ord}_t \cI : \cL (X) \rightarrow \NN \cup \{\infty\}$ given by
$\varphi \mapsto \min_{g} {\rm ord}_t g (\varphi)$,
where the minimum is taken over all $g$ in the stalk $\cI_{\pi_{0}
(\varphi)}$
of $\cI$ at $\pi_{0}
(\varphi)$. 
Let $\Omega^{1}_{X}$ be the sheaf of differentials on $X$ and let
$\Omega^{d}_{X}$ be the
$d$-th exterior power of 
$\Omega^{1}_{X}$. If $X$ is smooth and $\cF$ is a coherent sheaf
on $X$ together with a natural morphism
$\iota : \cF \rightarrow \Omega^{d}_{X}$, we denote by $\cI (\cF)$ the sheaf of ideals on $X$
which is locally generated by functions $\iota (\omega) / dx$ with $\omega$ a
local section of $\cF$ and $dx$ a local generator of $\Omega^{d}_{X}$.
For $n \in \NN \cup \{\infty\}$, we shall denote by 
${\rm ord}_t \cF = n$ the definable subassignement
$$
K \longmapsto \Bigl\{x \in X (K[[t]]) \Bigm|
\ord_{t} \cF \, (x) = n
\Bigr\}.
$$
The subassignements ${\rm ord}_t \cF = n$, $n \in \NN \cup \{\infty\}$, form a definable
partition
of the functor $h_{\cL (X)}$.

\begin{theorem}\label{change}Let $X$ and $Y$ be irreducible
algebraic varieties 
over $k$ of dimension $m$. Assume $Y$ is smooth. Let $p : Y \rightarrow 
X$ be a proper birational morphism.
Let $h_{n}$ be a definable partition of a definable subassignement $h$
of $h_{\cL (X)}$
with parameter $n \in \NN$. Consider
the definable subassignement $k_{n}$
of $h_{\cL (Y)}$ defined by 
$$
k_{n} :=\bigcup_{i + j = n} p^{-1} (h_{i})
\cap ({\rm ord}_t p^{\ast} \Omega^{d}_{X} = j).
$$
The $k_{n}$'s, $n \in \NN \cup \{\infty\}$,
form a definable partition of 
$p^{-1} (h)$ and
$$
\int \LL^{- n} h_{n} d\nu_{X} =
\int \LL^{- n} k_{n} d\nu_{Y}.$$
Furthermore, if $h \cap h_{\cL (p (E))} = \emptyset$,
with $E$ the exceptional locus of $p$, and $h$ and the $h_{n}$'s are
weakly stable (hence stable by Remark \ref{4.1}), then $p^{-1} (h)$ and the
$k_{n}$'s are stable and the above formula still holds when $\nu$ is
replaced by $\tilde \nu$.
\end{theorem}

\begin{proof}Follows from Lemma 3.4 in \cite{Arcs} similarly
as the proof of Lemma 3.3 in \cite{Arcs}. 
\end{proof}

\section{Rationality results}\label{rs}
In this section we prove
analogues of rationality results in \cite{Arcs}.
\subsection{}We consider the ring
$\widehat K_{0}^{v} ({\rm Mot}_{k, \bar \QQ})_{\QQ} [[T]]$
of power series in the 
variable $T = (T_{1}, \ldots, T_{r})$ with coefficients in
$\widehat K_{0}^{v} ({\rm Mot}_{k, \bar \QQ})_{\QQ}$.
We
denote by
$\overline K_{0}^{v} ({\rm Mot}_{k, \bar \QQ})_{{\rm loc}, \QQ} [[T]]_{{\rm rat}}$
the subring of
$\widehat K_{0}^{v} ({\rm Mot}_{k, \bar \QQ})_{\QQ} [[T]]$
generated by 
$\overline K_{0}^{v} ({\rm Mot}_{k, \bar \QQ})_{{\rm loc}, \QQ} [T]$,
$(\LL^i - 1)^{-1}$ and $(1 - \LL^{-a} \, T^b)^{-1}$,
with $i \in \NN \setminus \{0\}$, $a \in \NN$, and
$b \in \NN^r \setminus \{0\}$.

\begin{theorem}\label{rat1}Let $X$ be an algebraic variety
over $k$ of  dimension $d$. Let $h_{n}$ be a
definable subassignement of $h_{\cL (X)}$ depending on the parameter 
$n \in \NN^{r}$, and let
$h_{n, i}$ be a definable partition of the $h_{n}$'s depending
on the parameter $(n, i)
\in \NN^r \times \NN$.
Then the power series
\begin{equation}
\sum_{n \in \NN^r} \int \LL^{-i} \, h_{n, i} \, d \nu \, T^n
\end{equation}
in the variable $T = (T_{1}, \ldots, T_{r})$
belongs to $\overline K_{0}^{v} ({\rm Mot}_{k, \bar \QQ})_{{\rm loc}, \QQ} [[T]]_{{\rm
rat}}$.
\end{theorem}

\begin{cor}\label{deno}For every definable subassignement $h$ of 
$h_{\cL (X)}$, the measure
$\nu (h)$ belongs to the subring
$\overline K_{0}^{v} ({\rm Mot}_{k, \bar \QQ})_{{\rm loc}, \QQ}
[((\LL^i - 1)^{-1})_{i \geq 1}]$ of
$\widehat K_{0}^{v} ({\rm Mot}_{k, \bar \QQ})_{\QQ}$.
\end{cor}

\subsection{}Let us
denote by
$K_{0}^{v} ({\rm Mot}_{k, \bar \QQ})_{{\rm loc}, \QQ} [[T]]_{{\rm rat}}$
the subring 
of $ K_{0}^{v} ({\rm Mot}_{k, \bar \QQ})_{{\rm loc}, \QQ} [[T]]$
generated by 
$ K_{0}^{v} ({\rm Mot}_{k, \bar \QQ})_{{\rm loc}, \QQ} [T]$ and the series
$(1 - \LL^{-a} T ^b)^{-1}$,
with $a \in \NN$ and
$b \in \NN^r \setminus \{0\}$.

\begin{theorem}\label{rat2}Let $X$ be an algebraic variety
over $k$ of  dimension $d$. Let $h_{n}$ be a
definable subassignement of  $h_{\cL (X)}$ depending on the parameter 
$n \in \NN^{r}$, and let
$h_{n, i}$ be a definable partition of the $h_{n}$'s depending
on the parameter $(n, i)
\in \NN^r \times \NN$.
Assume $h_{n} \cap h_{\cL (X_{0})} = \emptyset$, for every
$n \in \NN^{r}$,
with $X_{0}$ the union of the singular locus of $X$ and its irreducible
components of dimension $<d$ and that the
subassignements $h_{n}$ and 
$h_{n, i}$ are all weakly  stable (hence stable by Remark \ref{4.1}).
Then the power series
\begin{equation}\label{ser}
\sum_{n \in \NN^r}  \int \LL^{-i} \, h_{n, i} \, d  \tilde \nu \, T^n
\end{equation}
in the variable $T = (T_{1}, \ldots, T_{r})$
belongs to $K_{0}^{v} ({\rm Mot}_{k, \bar \QQ})_{{\rm loc}, \QQ} [[T]]_{{\rm rat}}$.
\end{theorem}

\subsection{}For the proof of Theorem \ref{rat2} one needs the following
technical
lemma on bounded representations.

We shall say that a definable subassignement $h_{\ell}$ depending on parameters
$\ell \in \NN^{n}$, of $h_{\cL (X)}$ 
has a bounded representation
if there exists a covering of $X$ by
affine Zariski open sets $X_{i}$ embedded
in $\AA^{m_{i}}_{k}$ such that
$h_{\ell} \cap h_{\cL (X_{i})}$ is associated to a special formula
on $k [[t]]^{m_{i}}$
which is 
obtained by repeated application of
conjunction and negation from formulas of the form
(\ref{aaa}),  (\ref{bbb}) and  (\ref{ccc})
with
$\ord_{t}f_{i}$ bounded on $h_{\ell } \cap
h_{\cL (X_{i})}$ for each fixed $\ell$.

Clearly, if the family $h_{\ell}$ has a bounded representation then 
each 
$h_{\ell}$ is weakly stable.

\begin{lem}\label{2.8}Let $X$ be a quasi-projective algebraic variety over 
$k$ and let $h_{\ell}$, $\ell \in \NN^{n}$, be a definable
subassignement of $h_{\cL (X)}$
depending on parameters $\ell \in \NN^{n}$. Assume that 
$h_{\ell}$ is weakly stable for each $\ell$. Then the family 
$h_{\ell}$ is a finite boolean combination of definable
subassignements of $h_{\cL (X)}$
depending on parameters
which
have bounded representations.
\end{lem}

\begin{proof}The proof of Lemma 2.8 of \cite{Arcs} may be directly
adapted to carry over to the present situation.
\end{proof}

\subsection{}\label{pro}The proofs of Theorem \ref{rat1} and \ref{rat2}
are quite similar to the ones of Theorem 5.1 and
Theorem 5.1$'$ in \cite{Arcs}. We give details which will be 
used in section \ref{interpol}.

\begin{proof}[Proof of Theorem \ref{rat1} and \ref{rat2}]Let us first prove
Theorem \ref{rat2}. By Theorem \ref{change} one may use a
resolution  of singularities and assume that $X$ is smooth. Also we may
assume
$X$ is affine. By Pas's Theorem, Corollary \ref{prelq}, and Lemma \ref{2.4},
there exists a Presburger function $\theta : \NN^{r} \rightarrow \NN$
such that the series (\ref{ser}) is equal to
\begin{equation}
\sum_{i \leq \theta (n), n \in \NN^{r}} \tilde \nu (h_{n, i}) \, \LL^{-i} \,
T^{n}.
\end{equation}
By Lemma \ref{2.8} one reduces as in \cite{Arcs} to the case where the
family
$(h_{n, i})_{(n, i)}$ has bounded representation. Furthermore,
one may assume $X = X_{i}$ in the bounded representation of the family 
$(h_{n, i})_{(n, i)}$ by special formulas
obtained by repeated application of
conjunction and negation from formulas of the form
(\ref{aaa}),  (\ref{bbb}) and  (\ref{ccc}).
We denote by $F$ the product of all the polynomials
$f_{i}$, assumed to be non zero, occuring in these formulas
of the form
(\ref{aaa}),  (\ref{bbb}) and  (\ref{ccc}) and
we consider an embedded resolution of singularities
$\gamma : Y \rightarrow X$ of the locus 
of $F= 0$ in $X$, with exceptional locus contained in $\gamma^{-1} (F^{- 1} (0))$.
The variety $Y$ admits a covering by affine open subsets $U$ on which
there exist regular functions
$z_1, \ldots, z_d$ inducing an {\'e}tale map
$U \rightarrow \AA^d_k$ such that on $U$
each $f_i \circ \gamma$ is a monomial in
$z_1, \ldots, z_d$ multiplied by a regular function with no zeros on 
$U$. One may assume furthermore that the variables
$z_{i}$ appearing in at least one of these monomials
are exactly $z_{1}$, $z_{2}$, \dots, $z_{d_{0}}$.

Now for $w$ a definable
subassignement of $h_{(\AA^{1}_{k} \setminus
\{0\})^{d_{0}} \times U}$ and $\ell_{1}, \ldots, \ell_{d_{0}}$ in $\NN$,
we denote by
$h_{w, \ell_{1}, \ldots, \ell_{d_{0}}}$ the definable
subassignement of $h_{\cL (U)}$ defined by
\begin{multline}
K \longmapsto\\
\bigl\{ x \in \cL (U) (K) \bigm|
\ord_{t} z_{i} (x) = \ell_{i}, 1 \leq i \leq d_{0}
\quad
\textrm{and}  \quad
((\ac (z_{i} (x))), \pi_{0} (x)) \in w (K)
\bigr\}.
\end{multline}

It now follows from Theorem \ref{change} and the fact that 
$\ord_{t} F$ is bounded on $h_{n, i}$, that,
uniformly in $n$, $i$,
$\tilde \nu (h_{n, i})$ is a finite 
$\ZZ$-linear combination of
terms of the form
\begin{equation}\label{vol}
\sum_{{\ell_1, \ldots, \ell_{d_{0}} \in \NN}\atop{\theta (\ell_1,
\ldots, \ell_{d_{0}},n, i)}} 
\LL^{- \beta (\ell_1,
\ldots, \ell_{d_{0}})} \, \tilde \nu
(h_{w, \ell_1, \ldots, \ell_{d_{0}}}),
\end{equation}
where
$\theta (\ell_1,
\ldots, \ell_{d_{0}},n, i)$ is a condition defining
a Presburger subset of
$\ZZ^{d_{0} + r + 1}$,
$\beta$
is a linear form with coefficients in $\NN$, and
$w$ is a definable
subassignement of $h_{(\AA^{1}_{k} \setminus
\{0\})^{d_{0}} \times U}$, for $U$ as above.
Since 
$\ord_{t} F$ is bounded on $h_{n, i}$, the sum in (\ref{vol}) is finite.

Let us denote by $\bar w$ the definable
subassignement of $h_{(\AA^1_{k} \setminus 
\{0\})^{d_{0}} \times U}$
defined as
$\bar w = w \cap h_{U'}$ with $U'$ the subvariety of
$(\AA^1_{k} \setminus 
\{0\})^{d_{0}} \times U$ defined as the locus of 
points
$(w_{1}, \ldots,w_{d_{0}}, y)$
in $(\AA^1_{k} \setminus 
\{0\})^{d_{0}} \times U$
such that $z_{i} (y) = 0$
when $\ell_{i} > 0$ and $z_{i} (y) = w_{i}$
when $\ell_{i} = 0$. But, by 
Lemma 4.1 of \cite{Arcs}, with $n = e = 0$,  we have
\begin{equation}
\tilde \nu (h_{w, \ell_1, \ldots, \ell_{d_{0}}}) =
\chi_{c} (\bar w) \, \LL^{-(\sum_{i = 1}^{d_{0}} \ell_{i}) -d},
\end{equation} hence we can rewrite,
uniformly in $n$, $i$,
$\tilde \nu (h_{n, i})$ as a finite 
$\ZZ$-linear combination of
terms of the form
\begin{equation}\label{vol'}\LL^{-d}
\sum_{{\ell_1, \ldots, \ell_{d_{0}} \in \NN}\atop{\theta (\ell_1,
\ldots, \ell_{d_{0}},n, i)}} 
\LL^{- \beta (\ell_1,
\ldots, \ell_{d_{0}})} \, \chi_{c} (\bar w),
\end{equation}
with $\theta$, $\beta$ and $\bar w$ as above.

As in \cite{Arcs}, one may now
conclude the proof by using Lemma 5.2 and Lemma 5.3 of 
\cite{Arcs}.

\medskip

The proof of 
Theorem \ref{rat1} is similar except for the fact we have to replace
$\tilde \nu$ by $\nu$ and the finite sums are replaced by infinite sums
which converge in 
$\widehat K_{0}^{v} ({\rm Mot}_{k, \bar \QQ})_{\QQ} [[T]]$.
In particular,
$\nu (h_{n, i})$ is still a finite 
$\ZZ$-linear combination of terms of the form (\ref{vol'}), but the
number of terms in the series (\ref{vol'}) may now be infinite.
\end{proof}

\section{Arithmetic motivic integration specializes to $p$-adic integration}\label{interpol}

\subsection{}Let again $U = {\rm Spec} \, R$
denote an affine scheme of finite 
type over $\ZZ$, which will assumed to be integral and normal,
with fraction field
$k$ of characteristic zero.

\begin{lem}The morphism 
$$
{\rm \acute{E}t}_{\ell} : K_{0}^{v} ({\rm Mot}_{k, \bar \QQ})_{{\rm
loc}, \QQ}
\longrightarrow K_{0} (\bar \QQ_{\ell}, G_{k}) \otimes
\QQ
$$
induced by {\'e}tale $\ell$-adic 
realization factorizes
through a morphism
$$
{\rm \acute{E}t}_{\ell} :
\overline K_{0}^{v} ({\rm Mot}_{k, \bar \QQ})_{{\rm loc}, \QQ}
\longrightarrow K_{0} (\bar \QQ_{\ell}, G_{k}) \otimes
\QQ.
$$
In particular, for $x$ in $U$,
the morphism
$$
{\rm Tr} \, {\rm Frob}_{x} : K_{0}^{v} ({\rm Mot}_{k, \bar \QQ})_{{\rm
loc}, \QQ}
\longrightarrow \bar \QQ_{\ell} 
$$
factorizes
through a morphism
$$
{\rm Tr} \, {\rm Frob}_{x} :
\overline K_{0}^{v} ({\rm Mot}_{k, \bar \QQ})_{{\rm loc}, \QQ}
\longrightarrow \bar \QQ_{\ell}.
$$
\end{lem}

\begin{proof}For $\alpha$ in
$K_{0}^{v} ({\rm Mot}_{k, \bar \QQ})_{\rm loc}$,
denote by $P_{\alpha, x} (T)$ the ``characteristic polynomial''
${\rm det} (1 - {\rm Frob}_{x}T)$ of
${\rm Frob}_{x}$
on ${\rm \acute{E}t}_{\ell} (\alpha)$
(since we are dealing with virtual representations $P_{\alpha, x} (T)$
is a rational
function). Assume
$\alpha$ is in
$\cap_{m} F^{m} K_{0}^{v} ({\rm Mot}_{k, \bar \QQ})_{\rm loc}$ for
every $m$ in $\ZZ$.
By the part of the
Weil conjectures proven by Deligne \cite{Del}, we have then
$P_{\alpha, x} (T) = 1$ for $x$ in a set of Dirichlet density 1.
Since a virtual $\ell$-adic
representation of $G_k$ is determined by the
corresponding characteristic polynomials, we deduce from
Chebotarev's Theorem that ${\rm \acute{E}t}_{\ell} (\alpha) = 1$.
The result follows.
\end{proof}

We may extend  the morphism
$$
{\rm Tr} \, {\rm Frob}_{x} :
\overline K_{0}^{v} ({\rm Mot}_{k, \bar \QQ})_{{\rm loc}, \QQ}
\longrightarrow \bar \QQ_{\ell}
$$
to
a morphism
$$
{\rm Tr} \, {\rm Frob}_{x} : \overline K_{0}^{v} ({\rm Mot}_{k, \bar
\QQ})_{{\rm loc}, \QQ}
[((\LL^i - 1)^{-1})_{i \geq 1}]
\longrightarrow  \bar \QQ_{\ell},
$$
by sending
$(\LL^i - 1)^{-1}$ to $(q_{x}^i - 1)^{-1} $,
with $q_{x}$ the cardinality of $\FF_{x}$.

\begin{remark}Let  $\alpha$ be in one of the rings
$K_{0}^{v} ({\rm Mot}_{k, \bar \QQ})_{{\rm loc}, \QQ}$, 
$\overline K_{0}^{v} ({\rm Mot}_{k, \bar \QQ})_{{\rm loc}, \QQ}$, or 
$\overline K_{0}^{v} ({\rm Mot}_{k, \bar
\QQ})_{{\rm loc}, \QQ}
[((\LL^i - 1)^{-1})_{i \geq 1}]$.
For $x$ in a dense open subset of $U$, 
${\rm Tr} \, {\rm Frob}_{x} (\alpha)$ is in fact a rational number,
since it is a $\QQ$-linear
combination of number of rational
points of varieties over the residual field $\FF_{x}$.
In particular, it follows from 
Corollary \ref{deno} 
that, for every definable subassignement $h$ of 
$h_{\cL (X)}$, ${\rm Tr} \, {\rm Frob}_{x} (\nu (h))$ is a
rational number for $x$ in a dense open subset of $U$.
\end{remark}

Similarly,
the morphism 
${\rm Tr} \, {\rm Frob}_{x}$ may be naturally
extended 
to a morphism
$$
{\rm Tr} \, {\rm Frob}_{x} :
\overline K_{0}^{v} ({\rm Mot}_{\QQ, \bar \QQ})_{{\rm loc}, \QQ}
[[T]]_{{\rm rat}}
\longrightarrow \bar \QQ_{\ell} (T).
$$

\subsection{}Let $K$ be the field of fractions of a complete
discrete valuation ring
$\cO_{K}$ with finite residue field $\FF_{K}$.
We set $q = {\rm card} \, \FF_{K}$ and we fix an uniformizing parameter
$\pi$. We normalize the valuation $\ord$ on $K$ by $\ord (\pi) = 1$,
and we define, for $x \not = 0$ in $K$, $\ac x = \pi^{- \ord (x)} x$
and $|x| = q^{- \ord (x)}$, and $\ac 0  = |0| = 0$.
We obtain in this way a structure for the language
$\cL_{\rm Pas}$. 
By a definable subset of $\cO_{K}^{n} = \AA^{n} (\cO_{K})$,
we shall mean a subset defined by
an $\cL_{\rm Pas}$-formula with coefficients in $K$ in the valued
field sort and coefficients in $\FF_{K}$ in the residue field sort, $n$ free variables
running over $K$ and no other free variables. 
More generally, for $\cX$ an algebraic variety over $\cO_{K}$,
one defines definable subsets of $\cX (\cO_{K})$ by taking covers by affine open
subvarieties. Similarly, one defines definable subsets of $\cX (\cO_{K})$
depending on parameters in $\ZZ^{r}$, $\NN^{r}$, or $\NN \cup
\{\infty\}$.

Let $d$ be the Krull dimension of $\cX \otimes K$.
There is a natural $d$-dimensional measure on $\cX (\cO_{K})$,
cf. \cite{Serre}, \cite{Oesterle}, which we shall denote by $\nu_{K}$, for which
all definable subsets of $\cX (\cO_{K})$ are measurable, when $K$ is of
characteristic zero, cf.
\cite{cell}, \cite{dpp}, \cite{Veys}. This measure is defined by requiring
that the fibers of the reduction map modulo $\pi$ have measure $q^{-d}$.

\subsection{}We assume now $k$ is a finite
extension of $\QQ$ with ring of integers $\cO$ and $R = \cO
[\frac{1}{N}]$,
for some non zero integer $N$ which is a multiple of the discriminant of
$k$.
For $x$ a closed point of $\spec R$,
we denote by $K_{x}$ the completion of the localization of $R$ at $x$,
by $\cO_{K_{x}}$ its ring of integers
and by $\FF_{x}$ the residue field at $x$.

Let $\varphi$ be a formula over $R [[t]]^{m}$. 
We denote by $Z (\varphi, \cO_{K_{x}})$ 
the subset of $\cO_{K_{x}}^{m}$
defined by the formula $\varphi$.
By Proposition \ref{axax}, if $\varphi'$ is another
formula over $R [[t]]^{m}$ defining the same definable
subassignement of $h_{\cL (\AA^{m}_{R})}$ as 
$\varphi$, we have 
$$Z (\varphi, \cO_{K_{x}}) = Z (\varphi', \cO_{K_{x}})$$
for almost all $x$. Now let $\cX$ be a variety over
$R$ and let $h$ be a definable subassignement of
$h_{\cL (\cX)}$. Take a covering of $\cX$ by affine open subvarieties $U_{i}$
together with 
formulas $\varphi_{i}$ 
defining $h$ on $U_{i}$. By the preceding observation,
for almost all closed point $x$ of $\spec R$,
the subsets $Z (\varphi_{i}, \cO_{K_{x}})$ of 
$U_{i} (\cO_{K_{x}})$ may be glued together to define
a subset of $\cX (\cO_{K_{x}})$, which will
denote by $h [\cO_{K_{x}}]$. This subset
$h [\cO_{K_{x}}]$ is well
defined only for almost all $x$.

\medskip

As the following result shows, $p$-adic integration
may be viewed as a  specialization of 
arithmetic motivic integration.

\begin{theorem}\label{compa}Let $k$ be a finite
extension of $\QQ$ with ring of integers $\cO$ and $R = \cO
[\frac{1}{N}]$,
for some non zero integer $N$.
Let $\cX$ be a variety over $R$ and let $h$ be a
definable subassignement of
$h_{\cL (\cX)}$. Then there exists a non zero multiple $N'$ of 
$N$, such that, for every closed point $x$ of
$\spec \cO[\frac{1}{N'}]$,
$${\rm Tr} \, {\rm Frob}_{x}
\Bigl(
\nu (h \otimes k)
\Bigr)
=
\nu_{K_{x}} \Bigl(h [\cO_{K_{x}}]\Bigr).$$
More generally, 
if $h_{n}$ is a
definable subassignement of $h_{\cL (\cX)}$ depending on the parameter 
$n \in \NN^{r}$, and if
$h_{n, i}$ is a definable partition of the $h_{n}$'s depending
on the parameter $(n, i)
\in \NN^r \times \NN$,
there exists a non zero multiple $N'$ of 
$N$, such that, for every closed point $x$ of
$\spec \cO[\frac{1}{N'}]$, the series
$\sum_{i \in \NN} q_{x}^{-i}
\nu_{K_{x}} (h_{n, i} [\cO_{K_{x}}])$ converges in $\RR$ to a
rational number, for every $n$ in $\NN^{r}$, and
the following equality of power series 
\begin{equation}
{\rm Tr} \, {\rm Frob}_{x}
\Bigl(\sum_{n \in \NN^r}  \int \LL^{-i} (h_{n, i}  \otimes k)
d \nu \, T^n
\Bigr)
= \sum_{n \in \NN^r}  \sum_{i \in \NN} q_{x}^{-i}
\nu_{K_{x}} \Bigl(h_{n, i} [\cO_{K_{x}}]  \Bigr) \, T^n
\end{equation}
holds in $\QQ (T)$.
\end{theorem}

\begin{proof}We shall prove directly the more general second statement.
The proof will proceed by comparison with the proof of
Theorems \ref{rat1} and \ref{rat2}. We shall make here no notational
difference between $h$ and $h \otimes k$ and between $h_{n, i}$
and $h_{n, i} \otimes k$.
We set $X = \cX \otimes k$. In loc. cit. we first used a resolution
of singularities of $X$ to reduce to the smooth case, and then considered
an embedded resolution of a divisor $F = 0$
on some affine open subsets of $X$.
By Theorem 2.4 of \cite{AMJ}, there exists a
non zero multiple $N'$ of 
$N$ such that
these resolutions extend over $\spec \cO[\frac{1}{N'}]$
to resolutions with good reduction mod $P_{x}$, for every closed point
$x$ in
$\spec \cO[\frac{1}{N'}]$, in the sense of \cite{AMJ}. Here 
$P_{x}$ denotes the maximal ideal at $x$.
Now, by the local calculations of $p$-adic integrals on 
resolutions with good reduction in \cite{AMJ},
we
deduce that, for every closed point
$x$ in
$\spec \cO[\frac{1}{N'}]$,
$\nu_{K_{x}} (h_{n, i} [\cO_{K_{x}}])$
is a finite $\ZZ$-linear combination
of terms of the form
\begin{equation}\label{vol''}q_{x}^{-d}
\sum_{{\ell_1, \ldots, \ell_{d_{0}} \in \NN}\atop{\theta (\ell_1,
\ldots, \ell_{d_{0}},n, i)}} 
q_{x}^{- \beta (\ell_1,
\ldots, \ell_{d_{0}})} \, 
{\rm card} \, \bar w (\FF_{x})
\end{equation}
where $\theta$, $\beta$ and $\bar w$ are the same as in (\ref{vol'}),
and furthermore the coefficients of the terms
(\ref{vol''}) in $\nu_{K_{x}} (h_{n, i} [\cO_{K_{x}}])$
are the same as the coefficients of the terms
(\ref{vol'}) in $\nu (h_{n, i})$.
Indeed, we may assume $X$ is affine and, by 
Proposition \ref{speci}, we may also assume, maybe after replacing $N$ by
some non zero multiple, that
the subassignements $h_{n,i}$
are defined by repeated application of
conjunction and negation from formulas of the form
\begin{gather}
 \ord f_1 (x_1, \ldots, x_m) \geq
\ord f_2 (x_1, \ldots, x_m) + L (n, i) \label{aaa2}\\
\ord f_1 (x_1, \ldots, x_m) \equiv
L (n, i) \mod b \label{bbb2}
\intertext{and}
\vartheta (\ac (f_1 (x_1, \ldots, x_m)),
\dots,
\ac (f_{m'} (x_1, \ldots, x_m))), \label{ccc2}
\end{gather}
where $f_i$ 
are regular functions, $L$ is a
polynomial with coefficients in $\ZZ$ and degree $\leq 1$,
$b$ is in  $\ZZ$, and  $\vartheta$ is an $\LL_{k}$-formula
in $m'$ free variables. 
Furthermore we may assume 
there exist regular functions
$z_1, \ldots, z_d$ on $X$ inducing an {\'e}tale map
$X \rightarrow \AA^d_k$ such that on $X$
each $f_i$ is a monomial in
$z_1, \ldots, z_d$ multiplied by a regular function with no zeros on 
$X$, and that the variables
$z_{i}$ appearing in at least one of these monomials
are exactly $z_{1}$, $z_{2}$, \dots, $z_{d_{0}}$.

With the notations from \ref{pro}, we deduce from
the local calculations in \ref{vol'} that
$\nu (h_{n, i})$ is a $\ZZ$-linear combination of terms of the form
\begin{equation}\label{pasc}\LL^{-d}
\sum_{{\ell_1, \ldots, \ell_{d_{0}} \in \NN}\atop{\theta (\ell_1,
\ldots, \ell_{d_{0}},n, i)}} 
\LL^{- \beta (\ell_1,
\ldots, \ell_{d_{0}})} \, \chi_{c} (\bar w),
\end{equation}
with 
$\theta (\ell_1,
\ldots, \ell_{d_{0}},n, i)$  a condition defining
a Presburger subset of
$\ZZ^{d_{0} + r + 1}$,
$\beta$
a linear form with coefficients in $\NN$, and
$\bar w$ a definable
subassignement of $h_{(\AA^{1}_{k} \setminus
\{0\})^{d_{0}} \times X}$. 
Here $\theta$ and $\beta$  depend only on  the
monomials appearing in the functions $f_{i}$, the coefficients of the
polynomials $L$ and the integers $b$, and $\bar w$ depends furthermore on
$\vartheta$.
In the good reduction case, local calculation of $p$-adic integrals
as performed in \cite{AMJ} just provides
the same expression for 
$\nu_{K_{x}} (h_{n, i} [\cO_{K_{x}}])$  as a 
$\ZZ$-linear combination of terms of the form
$$q^{-d}
\sum_{{\ell_1, \ldots, \ell_{d_{0}} \in \NN}\atop{\theta (\ell_1,
\ldots, \ell_{d_{0}},n, i)}} 
q^{- \beta (\ell_1,
\ldots, \ell_{d_{0}})} \, {\rm card} \,  \bar w (\FF_{x}),$$
with the same coefficients in $\ZZ$ as for $\nu (h_{n,i})$, and for each term
the same functions
$\theta$ and $\beta$ and the same definable
subassignement $\bar w$ as in (\ref{pasc}).
The result now follows, since, by
Proposition \ref{ppp},
$${\rm Tr} \, {\rm Frob}_{x} (\chi_{c} (\bar w)) = 
{\rm card} \, \bar w (\FF_{x}),$$
for every  closed point $x$  in
$\spec \cO[\frac{1}{N''}]$, for a suitable non zero multiple $N''$ of $N'$.
\end{proof}

\medskip
The following variant of Theorem \ref{compa} is proved in
the same way.

\begin{theorem}\label{compabis}Let $k$ be a field of characteristic zero
which is the field of fractions of a normal domain $R$
of finite type over $\ZZ$.
Let $\cX$ be a variety over $R$ and let $h$ be a
definable subassignement of
$h_{\cL (\cX)}$. Then there exists a non zero element $f$ of 
$R$, such that, for every closed point $x$ of
$\spec R_{f}$, $h (\FF_{x})$ is $\nu_{\FF_{x}[[t]]}$-measurable and
$${\rm Tr} \, {\rm Frob}_{x}
\Bigl(
\nu (h \otimes k)
\Bigr)
=
\nu_{\FF_{x}[[t]]} \Bigl(h (\FF_{x})\Bigr).$$
More generally, 
if $h_{n}$ is a
definable subassignement of $h_{\cL (\cX)}$ depending on the parameter 
$n \in \NN^{r}$, and if
$h_{n, i}$ is a definable partition of the $h_{n}$'s depending
on the parameter $(n, i)
\in \NN^r \times \NN$, then there exists a non zero element $f$ of 
$R$, such that, for every closed point $x$ of
$\spec R_{f}$, all the sets
$h_{n, i} (\FF_{x})$ are $\nu_{\FF_{x}[[t]]}$-measurable,
the series
$\sum_{i \in \NN} q_{x}^{-i}
\nu_{\FF_{x}[[t]]} (h_{n, i} [\cO_{K_{x}}])$ converges in $\RR$ to a
rational number, for every $n$ in $\NN^{r}$, and
the following equality of power series 
\begin{equation}
{\rm Tr} \, {\rm Frob}_{x}
\Bigl(\sum_{n \in \NN^r}  \int \LL^{-i} (h_{n, i}  \otimes k)
d \nu \, T^n
\Bigr)
= \sum_{n \in \NN^r}  \sum_{i \in \NN} q_{x}^{-i}
\nu_{\FF_{x}[[t]]} \Bigl(h_{n, i} (\FF_{x})  \Bigr) \, T^n
\end{equation}
holds in $\QQ (T)$.\qed 
\end{theorem}

\subsection{}Again, let $k$ be a finite
extension of $\QQ$ with ring of integers $\cO$.
Consider a polynomial $f (x) \in 
k [x_{1}, \ldots, x_{m}]$
and let
$\varphi$
be an $\cL_{{\rm Pas}}$-formula with 
coefficients in $k$
in the valued field sort and in the residue
field sort, free
variables $x_{1}, \ldots, x_{m}$ running
over the valued field sort and no other free
variables.
Let $x$ be a closed point of $\spec \cO$ and
set
$$W_{x} := \{y \in \cO_{K_{x}}^m \, \vert \, \varphi (y) \, \text{holds}\}.$$
We consider the 
$p$-adic integral
$$
I_{\varphi, f} (s, x) := \int_{W_{x}} \vert  f  (y) \vert^s_{x} \vert 
dy\vert_{x},
$$
where $\vert \, \, \vert_{x}$ and $\vert 
dy\vert_{x}$ denote the $p$-adic norm on $K_{x}$ and  volume form on
$K_{x}^{m}$, respectively.
By a result of Denef \cite{rat}, $I_{\varphi, f} (s, x)$
is a rational function of $q_{x}^{-s}$. Macintyre
\cite{Macintyre}, Pas \cite{P}
and Denef \cite{AMJ} proved that
the degrees of the numerator and denominator of this rational
function are bounded independently of
$x$.

\subsection{}We will now show that a much stronger
uniformity statement holds: the integrals $I_{\varphi, f} (s, x)$
may be interpolated in a canonical way
by a ``motivic rational function".

Let  $k$ be a field of characteristic zero. Let 
$f (x)$ be a polynomial in
$k [x_{1}, \ldots, x_{m}]$, or more generally a
definable function in the valued field variables and with values in 
the valued field
with coefficients in $k$
in the language $\cL_{\rm Pas}$,
and let 
$\varphi$
be an $\cL_{{\rm Pas}}$-formula with 
coefficients in $k$
in the valued field sort and in the residue
field sort, free
variables $x_{1}, \ldots, x_{m}$ running
over valued field variables
and no other free
variables.
Consider, for $n \in \NN \cup \{\infty\}$,
the definable subassignement $h_{\varphi, n}$ of
$h_{\cL (\AA_{k}^{m})}$ associated to
the formula $\varphi \wedge \ord f = n$. The 
$h_{\varphi, n}$'s form a
definable partition of the
definable subassignement of $h_{\cL (\AA_{k}^{m})}$
defined by $\varphi$.
Now set 
$$
I_{\varphi, f, {\rm mot}} (T) := \sum_{n \in \NN}
\nu (h_{\varphi, n}) T^{n}.
$$

\begin{theorem}\label{red}Let  $k$ be a finite
extension of $\QQ$ with ring of integers $\cO$. Let 
$f (x)$ be a polynomial in
$k [x_{1}, \ldots, x_{m}]$, or more generally a
definable function in the valued field variables and with values in 
the valued field
with coefficients in $k$ 
in the language $\cL_{\rm Pas}$,
and let 
$\varphi$
be an $\cL_{{\rm Pas}}$-formula with 
coefficients in $k$
in the valued field sort and in the residue
field sort, free
variables $x_{1}, \ldots, x_{m}$ running
over valued field variables
and no other free
variables.
Then the series
$I_{\varphi, f, {\rm mot}} (T)$ is canonically associated to $\varphi$
and $f$, it is a rational function belonging to
$\overline K_{0}^{v} ({\rm Mot}_{k, \bar \QQ})_{{\rm loc}, \QQ} [[T]]_{{\rm rat}}$,
and, for almost all closed point $x$ in $\spec \cO$,
the equality
$$
{\rm Tr} \, {\rm Frob}_{x} (I_{\varphi, f, {\rm mot}} (T))
=
I_{\varphi, f} (s, x),
$$ holds after setting
$T = q_{x}^{-s}$.
Furthermore, if $\varphi'$
is another  $\cL_{{\rm Pas}}$-formula with 
coefficients in $k$ in the valued field sort and in the residue
field sort, free
variables $x_{1}, \ldots, x_{m}$ running
over valued field variables
and no other free
variables, then
\begin{enumerate}
\item[(1)]$I_{\varphi \vee \varphi', f, {\rm mot}} (T)
=
I_{\varphi, f, {\rm mot}} (T) + 
I_{\varphi', f, {\rm mot}} (T) -
I_{\varphi \wedge \varphi', f, {\rm mot}} (T)$
\item[(2)]If 
$\varphi \approx \varphi'$, $I_{\varphi, f, {\rm mot}} (T) = 
I_{\varphi', f, {\rm mot}} (T)$.
\end{enumerate}
\end{theorem}

\begin{proof}The rationality of $I_{\varphi, f, {\rm mot}} (T)$
follows 
from Theorem \ref{rat1} and the fact that
it satisfies the required specialization property follows from
Theorem \ref{compa}. Assertions (1) and (2) just follow from \ref{ami}
(2) and (3).
\end{proof}

\section{Arithmetic Poincar{\'e} series versus geometric
Poincar{\'e} series}\label{versus}
\subsection{The geometric Poincar{\'e} series $P_{\rm geom}$}\label{geom}
Let $k$ be a field
of characteristic 0 and 
let  $X$ be a variety over $k$.
Let us remind the reader that
we denote by $\pi_{n}$ the  morphism of 
truncation of arcs,
$$
\pi_{n} : \cL (X) \longrightarrow \cL_{n} (X).
$$

In the paper \cite{Arcs}, we considered
the Poincar{\'e} series
$$
P_{{\rm geom}, X} (T)
:=
\sum_{n \in \NN} \, [\pi_{n} (\cL (X))] \, T^n,
$$ 
and, more generally, for $W$ closed in $X$, the series
$$
P_{{\rm geom}, X, W} (T)
:=
\sum_{n \in \NN} \, [\pi_{n} (\cL (X) \cap \pi_{0}^{-1} (W))] \, T^n.
$$
One of our main results was the following 
rationality statement:

\begin{theorem}The series $P_{{\rm geom}, X, W} (T)$,
considered as a power series over the ring $K_{0} ({\rm Sch}_k)_{\rm loc}$,
belongs to the subring of
$K_{0} ({\rm Sch}_k)_{\rm loc} [[T]]$ generated by
$K_{0} ({\rm Sch}_k)_{\rm loc} [T]$ and the series
$(1 - \LL^a \, T^b)^{-1}$ with $a \in \ZZ$ and $b$ in $\NN \setminus 
\{0\}$.
\end{theorem}

\subsection{The arithmetic Poincar{\'e} series $P_{\rm ar}$}\label{ar}
Let $X$ 
be a variety defined over a field $k$ of characteristic 0.
Let $W$ be a closed subvariety of $X$ defined over $k$. 
Now we set
$$
P_{{\rm ar}, X, W} (T) := 
\sum_{n \in \NN} \, \chi_{c} (\pi_{n} (h_{\cL (X)} \cap \pi_{0}^{-1} (h_{W}))) 
\, T^n
$$
in 
$K_{0}^{v} ({\rm Mot}_{k, \bar \QQ})_{\QQ} [[T]]$.

\begin{theorem}\label{arrat}Let $h$ be a definable subassignement of
$h_{\cL (X)}$. Then the series
$$
\sum_{n \in \NN} \, \chi_{c} (\pi_{n} (h)) 
\, T^n
$$
belongs to $K_{0}^{v} ({\rm Mot}_{k, \bar \QQ})_{{\rm loc}, \QQ} [[T]]_{{\rm
rat}}$.
In particular,
the series
$P_{{\rm ar}, X, W} (T)$ belongs to 
$K_{0}^{v} ({\rm Mot}_{k, \bar \QQ})_{{\rm loc}, \QQ} [[T]]_{{\rm rat}}$.
\end{theorem}

\begin{proof}One reduces to the case where $X$ is a closed subvariety
of $\AA^{m}_{k}$. Then the result follows from Theorem \ref{rat2},
since
\begin{equation*}
\chi_{c} (\pi_{n, X} (h) ) \, \LL^{- (n+ 1) m}
=
\tilde \nu_{\cL (\AA^{m}_{k})} (\pi_{n, \AA^{m}_{k}}^{-1} (\pi_{n, X} (h)))
\end{equation*}
and 
the subassignements $\pi_{n, \AA^{m}_{k}}^{-1} (\pi_{n, X} (h))$ are stable
definable subassignements of 
$h_{\cL (\AA^{m}_{k})}$ depending on the parameter 
$n \in \NN$. Here $\pi_{n, X}$ and
$\pi_{n, \AA^{m}_{k}}$ denote the truncation morphisms
associated to $X$ and $\AA^{m}_{k}$, respectively.
\end{proof}

Assume now $k$ is a finite
extension of $\QQ$ with ring of integers $\cO$.
Let $R = \cO[\frac{1}{N}]$ with $N$ a non zero integer.
Let $\cX$ be a variety over $R $ and $\cW$
be  a subvariety. We set $X := \cX \otimes k$ and $W := \cW \otimes k$
and denote by $d$ the dimension of $X$.
For $x$ a closed point of $\spec R$, and $n \in \NN$,
we consider the canonical morphisms
$$
\pi_{n} : \cX (\cO_{K_{x}}) \longrightarrow \cX (\cO_{K_{x}} /
\pi_{x}^{n + 1} \cO_{K_{x}}),
$$
where $\pi_{x}$ is an uniformizing parameter of $\cO_{K_{x}}$,
and we denote by $N_{n, x} (\cX, \cW)$ the cardinality of the finite set
$\pi_{n}(\cX (\cO_{K_{x}}) \cap \pi_{0}^{-1} (\cW (\FF_{x}))) $.
By \cite{rat} the Poincar{\'e} series
$$
P_{x, \cX, \cW} (T) := \sum_{n \in \NN} N_{n, x} (\cX, \cW)  \, T^{n}
$$
is the expansion of a rational function in $\QQ (T)$.

\begin{theorem}\label{reduc}Let $k$ be a finite
extension of $\QQ$ with ring of integers $\cO$ and $R = \cO
[\frac{1}{N}]$,
for some non zero integer $N$.
Let $\cX$ be a variety over $R$ and let
$\cW$ be a 
subvariety. Then there exists a non zero multiple $N'$ of 
$N$, such that, for every closed point $x$ of
$\spec \cO[\frac{1}{N'}]$,
$${\rm Tr} \, {\rm Frob}_{x}
\Bigl(P_{{\rm ar}, X, W} (T)
\Bigr)
= P_{x, \cX, \cW} (T)
$$
in $\QQ (T)$.
\end{theorem}

\begin{proof}We may assume
$\cX$ is a closed subvariety
of $\AA^{m}_{R}$ defined by 
$f_{i} = 0$, $1 \leq i \leq r$, and that $\cW$ is defined by
$g_{i} = 0$, $1 \leq i \leq s$,
$h_{i} \not= 0$, $1 \leq i \leq t$, with
$f_{i}, g_{i}, h_{i}$ in 
$R [x_{1}, \ldots, x_{m}]$. We may also assume that $N$ is a multiple of
the discriminant of $k$.

Consider the following
formula $\varphi$ in $\cL_{{\rm Pas}}$ in the $m + 1$
free variables $x = (x_{1}, \ldots, x_{m})$ and $w$
over the valued field sort,
\begin{multline*}
\bigwedge_{i} (\ord (g_{i} (x)) > 0)
\wedge
\bigwedge_{i} (\ord (h_{i} (x)) = 0)\\
\wedge
\exists (y_{1}, \ldots, y_{m})
\Bigl(\bigwedge_{i} (\ord (x_{i} - y_{i}) \geq \ord (w)) \wedge
\bigwedge_{i} (f_{i} (y) = 0)\Bigr).
\end{multline*}

For every closed point $x$ in $\spec R$, the equality 
$$P_{x, \cX, \cW} (q_{x}^{- m - 1} \,q_{x}^{-s})
=
\frac{q_{x}}{q_{x}-1} I_{\varphi, w} (s, x)
$$
holds, as an direct calculation shows, cf. Lemma 3.1 of \cite{rat}.
Similarly, the equality
$$
P_{{\rm ar}, X, W} (\LL^{- m - 1} \,T)
=
\frac{\LL}{\LL-1}
I_{\varphi, w, {\rm mot}} (T)
$$
follows from the very definitions,
hence the result is a direct consequence of
Theorem \ref{red}.
\end{proof}

\section{An example:  branches of plane curves}\label{ex}
\subsection{}Let $k$ be a field of characteristic zero. Consider
a formal branch of plane curve
$$
X :
\begin{cases}
x = w^m&\\
y = \sum_{j \geq m} a_{j}w^j,&
\end{cases}
$$
with coefficients
$a_{j}$ in $k$.
We define $e_{0} = m$, 
$\beta_{1} = {\rm inf} \{ j \vert a_{j} \not= 0 \, \text{and} \, e_{0}
\nmid j\}$,
$e_{1} = (e_{0}, \beta_{1})$, and, by induction,
$\beta_{i} = {\rm inf} \{ j \vert a_{j} \not= 0 \, \text{and} \, 
e_{i - 1}
\nmid j\}$,
$e_{i} = (e_{i - 1}, \beta_{i})$. The sequence of the integers
$e_{i}$ being 
strictly decreasing,
for some smallest integer $g$, $e_{g} = 1$. 
We also define $n_{i}$, for $1 \leq i \leq g$, by
$e_{i - 1} = n_{i} e_{i}$, and set $N_{i} := \prod_{1 \leq j \leq i} n_{j}$.
We set $\beta_{0} = m$,
$\beta_{g + 1} = \{\infty\}$ and $N_{0} = 1$.

\medskip
For any integer $r \geq 1$, we denote by $\mu (r)$
the group of $r$-th roots of unity in $\CC$.
The following well known
lemma is stated for convenience.

\begin{lem}\label{wk}Let $\zeta$ be an $m$-th root of unity.
Then
$
{\rm ord}_{w} (y(w) - y(\zeta w)) = \beta_{i}
$
if and only if
$\zeta \in \mu (e_{i - 1}) \setminus \mu (e_{i})$, for $i \geq 1$. \hfil \qed
\end{lem}

\subsection{}The definition of $\cL (X)$, $\cL_{n} (X)$ and $\pi_{n}$
is extended in a straightforward way to formal varieties to
set $\cX_{n} := \pi_{n} (\cL_{n} (X) \cap \pi_{0}^{-1} (0))$.
We define the
Poincar{\'e} series
$$P_{{\rm geom}, X, 0} (T) := \sum_{n \geq 0} \, [\cX_{n}] \, T^n.$$
When $X$ is a branch of an algebraic curve this definition coincides
with the one given in \ref{geom}.

To define  $P_{{\rm ar}}$ in this situation we consider
the formula
$\varphi_{n}$ in the free variables $x_{1}, \ldots, x_{n}$ and
$y_{1}, \dots,
y_{n}$ 
$$
\exists w \, \Bigl(
\ord_{t} (x - w^{m}) > n \wedge 
\ord_{t} (y - \sum_{m \leq j \leq n} a_{j}w^j) >n \Bigr),
$$
with $x = \sum_{1 \leq i \leq n} x_{i} t^{i}$,
$y = \sum_{1 \leq i \leq n} y_{i} t^{i}$ and 
$w = \sum_{1 \leq i \leq n} w_{i} t^{i}$,
and set 
$$
P_{{\rm ar}, X, 0} (T) = \sum_{n \geq 0} \chi_{c} (\varphi_{n}) T^{n}.
$$
When $X$ is a branch of an algebraic curve this definition coincides
with the one given in \ref{ar}.

\begin{prop}
Let $k$ be a field of characteristic zero. Consider
a formal branch of plane curve
$$
X :
\begin{cases}
x = w^m&\\
y = \sum_{j \geq m} a_{j}w^j,&
\end{cases}
$$
with coefficients $a_{j}$ in $k$.
Assume $k$ contains all $m$-th roots of unity.
Then
$$P_{{\rm geom}, X, 0} (T) = 
\frac{1}{1 - T} +
\frac{\LL - 1}{1 - \LL T}\frac{T^{m}}{1 - T^{m}}
$$
and
$$
P_{{\rm ar}, X, 0} (T) =
\frac{1}{1 - T} +
\frac{\LL - 1}{1 - \LL T}
\Bigl[\frac{1}{m}\frac{T^{m}}{1 - T^{m}}
+
\sum_{1 \leq i \leq g} \frac{N_{i} - N_{i -1}}{m}
\frac{\LL^{\beta_{i} - m}T^{\beta_{i}}}{1 - \LL^{\beta_{i} - m}T^{\beta_{i}}}
\Bigr].
$$
In particular, the  poles of the rational function
$P_{{\rm ar}, X, 0} (T)$ of the form $T =\LL^{\alpha}$, $\alpha$
in $\QQ
\setminus \{0, -1\}$,
are exactly
$T = \LL^{\frac{m}{\beta_{i}} - 1}$, $1 \leq i \leq g$.
\end{prop}

\begin{proof}For every integer $\ell$ such that
$0 < \ell \leq \frac{n}{m}$, we set
$\cX_{n, \ell}  = \cX_{n} \cap \{(x, y) \vert
{\ord}_{t} x = \ell m\}$,
with $x = \sum_{1 \leq i \leq n} x_{i} t^{i}$ and 
$y = \sum_{1 \leq i \leq n} y_{i} t^{i}$.
Similarly, we consider the formula
$\varphi_{n, \ell} = \varphi_{n} \, \wedge \, (\ord_{t} x = \ell m)$.
So we have
$
[\cX_{n}] = 1 + \sum_{0 < \ell \leq \frac{n}{m}}[\cX_{n,\ell}]
$
and $\chi_{c} (\varphi_{n}) = 1 + \sum_{0 < l \leq \frac{n}{m}}
\chi_{c} (\varphi_{n, \ell})$, and 
the result follows directly from the next lemma.
\end{proof}

\begin{lem}\label{last}Fix an integer $\ell \leq \frac{n}{m}$.
Let $i$ be the
unique integer $0 \leq i  \leq g$ such that
$\frac{n}{\beta_{i + 1}} < \ell \leq \frac{n}{\beta_{i}}$.
Then 
$$
[\cX_{n, \ell}] = (\LL - 1) \,\LL^{n - \ell m} \quad
\text{and}
\quad
\chi_{c} (\varphi_{n, \ell})
=
\frac{N_{i}}{m} \,
(\LL - 1) \,\LL^{n - \ell m}.
$$
\end{lem}

\begin{proof}We consider the variety
$W = \GG_{m, k} \times \AA^{n - \ell m}_{k}$
with coordinate $w_{\ell}$ on the first factor
and $(w_{\ell + 1}, \ldots, w_{n - \ell m + \ell})$ on the second
factor. There is a morphism
$h : W \rightarrow \AA^{n}_{k} \times \AA^{n}_{k}$
sending a point $(w_{\ell}, \ldots, w_{n - \ell m + \ell})$ to
the first $n$ coefficients of the series
$w^{m}$ and
$\sum_{j \ell \leq n} a_{j} w^{j}$. Let $Z$ be the image of $W$ by
$h$. It is a locally closed subvariety of 
$\AA^{n}_{k} \times \AA^{n}_{k}$. It follows from Lemma \ref{wk} that
$h : W \rightarrow Z$ is a Galois cover with group
$G = \mu (e_{i})$ acting by multiplication on $W$. By construction,
$\chi_{c} (\varphi_{n, \ell})$ is equal to 
$\chi_{c} (W, \delta)$, with $\delta$ the central function on $G$
which takes
value 1 at the identity and zero elsewhere.
But now remark that 
$\chi_{c} (W, \alpha) = 0$ for $\alpha$ a non trivial irreducible
character of $G$
and that $\chi_{c} (W, \alpha) = (\LL - 1) \,\LL^{n - \ell m}$ when
$\alpha$ is the trivial character
(cf. Lemma 1.4.3 of \cite{Motivic}). The second equality
follows, since, when $\delta$ is expressed
as a $\QQ$-linear combination of irreducible
characters, the coefficient of the trivial character is $e_{i}^{-1} = 
\frac{N_{i}}{m}$. For the first equality
it is enough to remark that $Z = \cX_{n,
\ell}$ and that $[W / \mu (e_{i})] = (\LL - 1) \, \LL^{n - \ell m}$. 
\end{proof}

\begin{remark}The calculation of the $p$-adic Poincar{\'e}
series for branches of plane curves was carried out
in \cite{Bollaerts}.
\end{remark}

\bibliographystyle{amsplain}

\end{document}